\documentclass[trsc,nonblindrev]{informs3} % current default for manuscript submission

\OneAndAHalfSpacedXI

\usepackage{natbib}
 \bibpunct[, ]{(}{)}{,}{a}{}{,}%
 \def\bibfont{\small}%
 \def\bibsep{\smallskipamount}%

\TheoremsNumberedThrough     % Preferred (Theorem 1, Lemma 1, Theorem 2)

\EquationsNumberedThrough    % Default: (1), (2), ...

\usepackage{longtable}
\usepackage{bbding}
\usepackage{amsmath}
\usepackage{amssymb}
\usepackage{multirow}
\usepackage{booktabs}
\usepackage{colortbl}
\usepackage{fancyhdr}
\usepackage{indentfirst}
\usepackage{amsmath, amssymb}
\usepackage{makecell}
\usepackage{multirow}
\usepackage{array, threeparttable}
\usepackage{cite}
\usepackage{url}
\usepackage{breakurl}
\usepackage{footnote}
\usepackage{threeparttable}
\usepackage{hyperref}
\usepackage{cleveref}
\usepackage{setspace}
\usepackage{lineno,hyperref}
\usepackage{amsmath,amssymb,ifthen,url,graphicx,color,array,theorem,cases}
\usepackage{xcolor}

\hypersetup{hypertex=true,
            colorlinks=true,
            linkcolor=blue,
            anchorcolor=blue,
            citecolor=blue}

\begin{document}

\RUNAUTHOR{Lin et al.}
\RUNTITLE{Sparse-Dense Flight Copy-Based Interactive Mechanism}

%%%%%%%%%%%%%%%%

\TITLE{Sparse-Dense Flight Copy-Based Interactive Mechanism for Airline Integrated Recovery}

\ARTICLEAUTHORS{

%\AUTHOR{Anonymous Author(s)}

\AUTHOR{Jiajin Lin}
\AFF{School of Information Engineering, Taizhou Vocational College of Science and Technology, Taizhou 318020, China, \\School of Mathematics, Nanjing University of Aeronautics and Astronautics, Nanjing 210016, China, \\\EMAIL{jiajin\_lin@outlook.com} \URL{}}
\AUTHOR{Jianlin Jiang, Yan Gu$^*$, Yuzhen Guo}
\AFF{School of Mathematics, Nanjing University of Aeronautics and Astronautics, Nanjing 210016, China, \\ \EMAIL{\{jiangjianlin@nuaa.edu.cn, guyanmath@nuaa.edu.cn, guoyuzhen@nuaa.edu.cn\}\\ $^{*}$ Corresponding author} \URL{}}
\AUTHOR{Cheng-Lung Wu}
\AFF{School of Aviation, University of New South Wales, Sydney, Kensington, NSW 2052, Australia, \\\EMAIL{c.l.wu@unsw.edu.au} \URL{}}

} % end of the block

\ABSTRACT{%
\textbf{Abstract.} Flight recovery, aircraft rerouting, and passenger reallocation are critical to airline disruption management. To better capture their interdependencies, we consider these recovery phases from an integrated perspective, incorporating cruise speed control to enhance recovery performance. Although flight copies are widely used to model flight rescheduling times in disruption modelling, the resulting integrated model is computationally challenging due to the large number of decision variables, especially when speed control is considered. To address this, we propose an innovative interactive mechanism based on a new sparse-dense flight copy approach. This mechanism alternately adjusts flight schedules and aircraft routes on a sparse flight copy network, and readjusts flight times and reallocates passengers on a dense flight copy network, with adjustments on one network utilizing feedback from the other. This significantly reduces problem complexity compared to the conventional flight copy approach. To operationalize this mechanism, we develop an integrated flight, aircraft, and passenger recovery model (IFAPRM) and design a customized Benders decomposition method, where the master problem addresses flight rescheduling with aircraft routing and the subproblem handles flight retiming with passenger reallocation. Column generation is embedded in both problems to avoid explicit enumeration of routes and itineraries. To accelerate computation, we introduce techniques based on favourable structural properties of IFAPRM, including an effective feasibility certificate, scale management, and valid inequalities. Experiments on real-world instances demonstrate that the sparse-dense flight copy-based interactive mechanism outperforms the conventional approach, reducing runtime by over 68\% for 121-flight instances and yielding high-quality solutions for 372-flight instances where the conventional approach fails to produce a solution. The acceleration techniques reduce runtime by approximately 27\%-47\%. In essence, the proposed mechanism, together with its modelling method, algorithm, and acceleration techniques, provides a comprehensive methodology for integrated rescheduling problems under complex operations, with broad applicability to transportation and logistics.}

\KEYWORDS{Airline integrated recovery; Sparse-dense flight copy approach; Interactive mechanism; Cruise speed control; Benders decomposition; Column generation}

%\linenumbers

\maketitle
\section{Introduction}\label{Introduction}
\noindent
Managing irregular flights is a common challenge in the airline industry. Once irregular flights disrupt the normal operations, the Airline Operations Control Center (AOCC) of an airline needs to take recovery actions to repair the disrupted schedules (\citeauthor{wu2018airlinecapacity}, \citeyear{wu2018airlinecapacity}). The efficiency of the recovery process significantly impacts not only direct costs but also customer satisfaction. According to \citeauthor{li2022} (\citeyear{li2022}), the direct loss is approximately \$35 billion globally, and the potential loss caused by airline disruptions rises to \$60 billion when accounting for lost passenger productivity and downstream economic effects in supporting industries such as hospitality, business services, and tourism.

To handle the complexity of practical operations, the recovery process is typically divided into four phases: flight recovery, aircraft rerouting, crew rescheduling, and passenger reallocation, which are often executed sequentially in practice. However, these phases are highly interdependent. Flight recovery directly affects airline revenue and strongly influences the other recovery phases. Decisions on aircraft rerouting, in turn, impact both crew rescheduling and passenger reallocation. Moreover, passenger reallocation plays a crucial role in passenger satisfaction for airline services. The above relations mean that these recovery phases form a tightly coupled system. Therefore, airline disruption management should be implemented from an integrated perspective.

Cruise speed control serves as an effective recovery option, capable of mitigating delay propagation across the network by compressing the flying times of flights (\citeauthor{akturk2014}, \citeyear{akturk2014}). Beyond reducing propagated delays, it enhances flexibility in departure and arrival times, creating more opportunities for aircraft swaps during rerouting. Consequently, the recovery plans become more adaptable than without cruise speed control. However, increasing cruise speeds to reduce delays leads to increased fuel consumption and consequent higher carbon emissions and costs. Therefore, the trade-off between reducing flight delays and additional fuel consumption should be considered when controlling cruise speeds in disruption recovery planning (\citeauthor{arikan2016}, \citeyear{arikan2016}; \citeauthor{marla2017}, \citeyear{marla2017}). Furthermore, considering cruise speed control in airline disruption management will significantly expand the number of decision variables in the resulting model, especially when incorporated into the commonly used flight copy approach where a flight copy represents an alternative schedule option for the original flight at a different departure time, which substantially increases the problem scale and computational burden.

\subsection{Flight Recovery, Aircraft Rerouting, and Passenger Reallocation}\label{Section1.1}
\noindent
Flight recovery, aircraft rerouting, and passenger reallocation are usually addressed in a sequential manner (\citeauthor{marla2017}, \citeyear{marla2017}). Specifically, flight recovery and aircraft rerouting are resolved first, followed by passenger reallocation. Such a sequential approach simplifies the entire recovery process, but it may lead to suboptimal recovery plans or infeasible subsequent recovery phases, because it does not fully consider the interdependencies among these phases.

To overcome this limitation, some researchers have focused on integrating flight recovery, aircraft rerouting, and passenger reallocation. \citeauthor{bratu2006} (\citeyear{bratu2006}) presented a pioneering study that simultaneously recovered disrupted flight schedules, aircraft routes, and passenger itineraries. The authors proposed two passenger-centric models to balance operating costs and passenger itinerary recovery costs. \citeauthor{sinclair2016} (\citeyear{sinclair2016}) formulated an integrated aircraft and passenger recovery problem as a mixed-integer linear programming model, and solved it via a two-stage column generation-based heuristic method. In the first stage, they selected the promising arcs on the networks for aircraft routes and passenger itineraries to limit the scale of the integrated recovery problem. Then, they fixed the arcs selected in the first stage and obtained the final recovery plan by column generation. \citeauthor{zhang2016} (\citeyear{zhang2016}) proposed a novel three-stage sequential math-heuristic framework: the first stage generated a preliminary flight schedule and aircraft routes; the second stage adjusted departure/arrival times to minimize estimated flight delay and passenger itinerary disruption costs; and the third stage handled passenger reallocation. The final recovery plan was obtained by iteratively solving the latter two stages. It is noted that this framework did not sufficiently account for the interdependencies between the first stage and the latter two stages. Inspired by the framework in \citeauthor{zhang2016} (\citeyear{zhang2016}), this paper establishes an interactive mechanism between different stages to capture their interdependencies, thereby improving the recovery performance.

In addition to the integration of flight recovery, aircraft rerouting, and passenger reallocation, various combinations of recovery phases have been extensively studied, with diverse solution methods developed for the airline integrated recovery problems (e.g., \citeauthor{petersen2012}, \citeyear{petersen2012}; \citeauthor{hu2021integrated}, \citeyear{hu2021integrated}; \citeauthor{xu2023distributionally}, \citeyear{xu2023distributionally}; \citeauthor{zhong2024}, \citeyear{zhong2024}; \citeauthor{wang2025}, \citeyear{wang2025}; \citeauthor{Jiang2026}, \citeyear{Jiang2026}). More details about this topic can be found in \citeauthor{hassan2021} (\citeyear{hassan2021}) and \citeauthor{su2021} (\citeyear{su2021}). More recently, \citeauthor{wu2025} (\citeyear{wu2025}) further summarized the literature on airline recovery published between 2021 and 2024.

As one of the first studies to computationally solve the fully integrated recovery problem (i.e., flight, aircraft, crew, and passenger recovery), \citeauthor{petersen2012} (\citeyear{petersen2012}) constructed several path-based models for different recovery phases and developed a Benders decomposition framework in which new aircraft routes and crew pairings were generated via column generation. They exploited the model's structure and innovatively proposed two effective certificates for detecting infeasibility and suboptimality of the Benders subproblem during the column generation process. Their infeasibility certificate is a necessary condition, so the column generation process continues until the condition of the infeasibility certificate is satisfied or the column generation process is completed. Motivated by \citeauthor{petersen2012} (\citeyear{petersen2012}), we propose a feasibility certificate instead of their infeasibility certificate to check the feasibility of Benders subproblem. In contrast to \citeauthor{petersen2012} (\citeyear{petersen2012}), our feasibility certificate is a necessary and sufficient condition, so the feasibility of the Benders subproblem can be easily determined by checking this certificate. Furthermore, this certificate does not require executing the column generation process for feasibility checks, which significantly reduces computational effort and improves efficiency.

\subsection{Cruise Speed Control and Flight Copy Approach}
\noindent
Cruise speed control is widely used in disruption management to enhance the flexibility and performance of recovery plans (\citeauthor{arikan2016}, \citeyear{arikan2016}; \citeauthor{marla2017}, \citeyear{marla2017}). The studies about the application of cruise speed control can be divided into two categories: continuous cruise speed control and discrete cruise speed control.

Continuous cruise speed control explicitly takes the cruise speed of a flight as a continuous decision variable that can be chosen in a given interval. A representative work by \citeauthor{akturk2014} (\citeyear{akturk2014}) focused on aircraft rerouting, employing a realistic fuel burn function to calculate the fuel consumption resulting from controlling cruise speeds. To address the resulting highly nonlinear function, the authors reformulated it as conic quadratic inequalities, transforming the intractable model into a conic quadratic mixed integer programming (CQMIP) model solvable by commercial solvers. Such a reformulation strategy has been extended to other fields in the airline industry, like robust schedule design (\citeauthor{duran2015}, \citeyear{duran2015}; \citeauthor{gurkan2016}, \citeyear{gurkan2016}) and airline integrated recovery (\citeauthor{arikan2016}, \citeyear{arikan2017}). While the reformulation strategy can address the realistic fuel burn function effectively, solving the CQMIP model remains challenging in practice, especially for the airline recovery problem on a large-scale airline network. To reduce the complexity of the nonlinear fuel consumption, \citeauthor{li2022} (\citeyear{li2022}) used a piecewise-linear function to approximate the nonlinear function.

As an approximation of continuous control, discrete cruise speed control selects the cruise speed of a flight from a discrete set. \citeauthor{marla2017} (\citeyear{marla2017}) implemented this by generating flight copies with varying flying times corresponding to different cruise speeds at each candidate departure time. \citeauthor{lee2020} (\citeyear{lee2020}) employed this flight copy approach to deal with the airline disruption management under airport operational uncertainty. \citeauthor{Zhang2024} (\citeyear{Zhang2024}) utilized discrete cruise speed control to address the operational aircraft maintenance routing problem (OAMRP). \citeauthor{Yuan2025} (\citeyear{Yuan2025}) also adopted this flight copy approach to solve the OAMRP with cruise speed control and further considered the balance of utilization rates across different resources.

The flight copy approach has been commonly used for modelling time-related decision-making in the airline industry, such as airline recovery (\citeauthor{maher2016}, \citeyear{maher2016}; \citeauthor{liang2018}, \citeyear{liang2018}; \citeauthor{Jiang2025}, \citeyear{Jiang2025}) and airline robust scheduling (\citeauthor{Froyland2014}, \citeyear{Froyland2014}; \citeauthor{xu2021}, \citeyear{xu2021}; \citeauthor{Schrotenboer2023}, \citeyear{Schrotenboer2023}). When applying this approach, both the solution quality and the runtime depend on the length of the copy interval for each flight. In general, shorter intervals offer more candidate flight times, but significantly increase decision variables and solution time, especially when the cruise speed control is incorporated. Such a trade-off becomes more pronounced as the problem's scale increases. To address this, several efficient solution methods have been proposed. \citeauthor{Huang2022} (\citeyear{Huang2022}) developed a copy evaluation method to avoid the generation of costly flight copies and designed a novel algorithm incorporating a copy generation and filtration process to control the number of involved copies. \citeauthor{Zhang2024} (\citeyear{Zhang2024}) reduced the excessive model size by deleting redundant copies with the same flight connection opportunities but higher costs.

\subsection{Research Scope and Aims}
\noindent
This paper focuses on the integration of flight recovery, aircraft rerouting, and passenger reallocation, which is inspired by two key studies: the three-stage math-heuristic recovery framework developed by \citeauthor{zhang2016} (\citeyear{zhang2016}) and the Benders decomposition framework proposed by \citeauthor{petersen2012} (\citeyear{petersen2012}). As noted in Section \ref{Section1.1}, \citeauthor{zhang2016} (\citeyear{zhang2016}) generated a preliminary flight schedule and aircraft routes in the first stage, then adjusted departure/arrival times and reallocated passengers alternately in the subsequent two stages via an alternating minimization method. Because this framework did not readjust the first-stage decision according to the outcomes of the subsequent stages, the final chosen flight copies may be unsatisfactory. In contrast, \citeauthor{petersen2012} (\citeyear{petersen2012}) proposed a Benders decomposition framework to iteratively select flight copies across all recovery phases. However, when the problem scale increases, a large number of decision variables corresponding to flight copies would be introduced, and then selecting flight copies by the Benders decomposition framework becomes challenging.

To address the limitations of the frameworks proposed by \citeauthor{zhang2016} (\citeyear{zhang2016}) and \citeauthor{petersen2012} (\citeyear{petersen2012}), this paper proposes an innovative sparse-dense flight copy-based interactive mechanism. First, building on the flight copy concept in \citeauthor{zhang2016} (\citeyear{zhang2016}), we introduce a new sparse-dense flight copy approach, which categorizes flight copies into two types according to the length of the copy interval: sparse flight copies (with long intervals) and dense flight copies (with short intervals), as illustrated by ``Level 1'' and ``Level 2'' in Figure \ref{sparse-denseFlightCopy}. Then, drawing inspiration from the Benders decomposition framework of \citeauthor{petersen2012} (\citeyear{petersen2012}), we develop an interactive mechanism that alternately selects sparse and dense flight copies within an integrated recovery model. This design captures interdependencies between recovery phases through interactive feedback, thereby enhancing recovery performance. Compared with the conventional flight copy approach, the proposed sparse-dense interactive mechanism demonstrates superior effectiveness.

\subsection{Contributions}
\noindent
The main contributions of this work are summarized as follows. First, we introduce a novel sparse-dense flight copy approach for the integrated flight, aircraft, and passenger recovery problem with cruise speed control. This integrated problem is then solved via a sparse-dense flight copy-based interactive mechanism. Under this mechanism, flight rescheduling and aircraft rerouting are performed on a sparse flight copy network, while flight retiming and passenger reallocation are conducted on a dense flight copy network. The final recovery plan is derived by alternately selecting sparse flight copies and dense flight copies from an integrated perspective. Essentially, the interactive mechanism establishes a bidirectional relationship between the two networks by adjusting the decision on one network according to the feedback information from the other. Consequently, the resulting recovery plan generally outperforms those produced by the sequential approach or heuristic methods, while the difficulty of finding an integrated solution is greatly reduced.

Second, to operationalize the interactive mechanism numerically, we develop a sparse-dense flight copy-based integrated recovery model, namely the integrated flight, aircraft, and passenger recovery model (IFAPRM), and propose a customized Benders decomposition (CBD) method to solve it. Notably, the process of alternately selecting sparse and dense flight copies aligns naturally with the Benders decomposition strategy of alternately solving Benders master problem and Benders subproblem. Besides, instead of simply addressing the two problems by an alternating minimization method, the CBD method seeks the optimal solution by adjusting the selection of sparse flight copies based on the decisions from the dense flight copy network, ensuring the entire recovery process is guided by an interactive perspective. Furthermore, the decomposition strategy significantly reduces problem complexity by dividing IFAPRM into two manageable problems.

Third, by exploiting the structural properties of the integrated recovery problem, we introduce several acceleration techniques to improve the efficiency of the CBD method. An effective feasibility certificate is first developed to reduce runtime in checking the feasibility of Benders subproblems. This certificate is a necessary and sufficient condition, allowing the feasibility of the Benders subproblem to be easily determined. Moreover, because this certificate operates independently of the column generation process, it eliminates the need for column generation in feasibility checks, thereby significantly improving efficiency. Additionally, based on known information during the iterative process of the CBD method, we remove redundant linking constraints when solving the Benders subproblem. By leveraging the specific structure of IFAPRM, we also develop stronger no-good cuts and stronger Laporte \& Louveaux cuts to tighten the formulation of the Benders master problem. Computational experiments in Section \ref{Section Numerical experiments} validate the effectiveness of these acceleration techniques, reducing the runtime by approximately 27\%-47\%.

Fourth, the sparse-dense flight copy-based interactive mechanism exhibits strong generality. The proposed mechanism, along with its corresponding modelling method, algorithm, and acceleration techniques, can address not only the integrated recovery problem but also serve as a comprehensive methodology and a general support framework for a wide range of decision-making problems with schedule adjustments. Moreover, within this methodology, decisions across different flight copy networks can be made using customized methods (e.g., heuristic algorithms for complex network problems and exact algorithms for small-scale constrained problems) tailored to the specific structure. This flexibility further validates the mechanism's broad applicability.

\subsection{Paper Outline}
\noindent
The remainder of this paper is organized as follows. Section \ref{Section Interactive mechanism} provides the relevant definitions for describing the integrated recovery problem and presents the sparse-dense flight copy-based interactive mechanism. Section \ref{Section Mathematical model} develops a mathematical model for the integrated recovery problem based on the interactive mechanism. Section \ref{Section Solution methodology} proposes a customized Benders decomposition method to solve the model. Section \ref{Section Acceleration techniques} introduces several acceleration techniques and presents the accelerated CBD method. Section \ref{Section Numerical experiments} evaluates the proposed methodology through extensive numerical experiments. The final section concludes this work and discusses future research directions. The proofs of properties are given in Appendix A.

\section{Interactive Mechanism}\label{Section Interactive mechanism}
\noindent
In this section, we provide the relevant definitions for the integrated flight, aircraft, and passenger recovery problem and then present the sparse-dense flight copy-based interactive mechanism.

\subsection{Definitions}\label{Section Definitions}
\begin{figure}[!htbp]
   \centering
   \begin{minipage}{1\linewidth}
    \centering
    \includegraphics[width=0.8\textwidth]{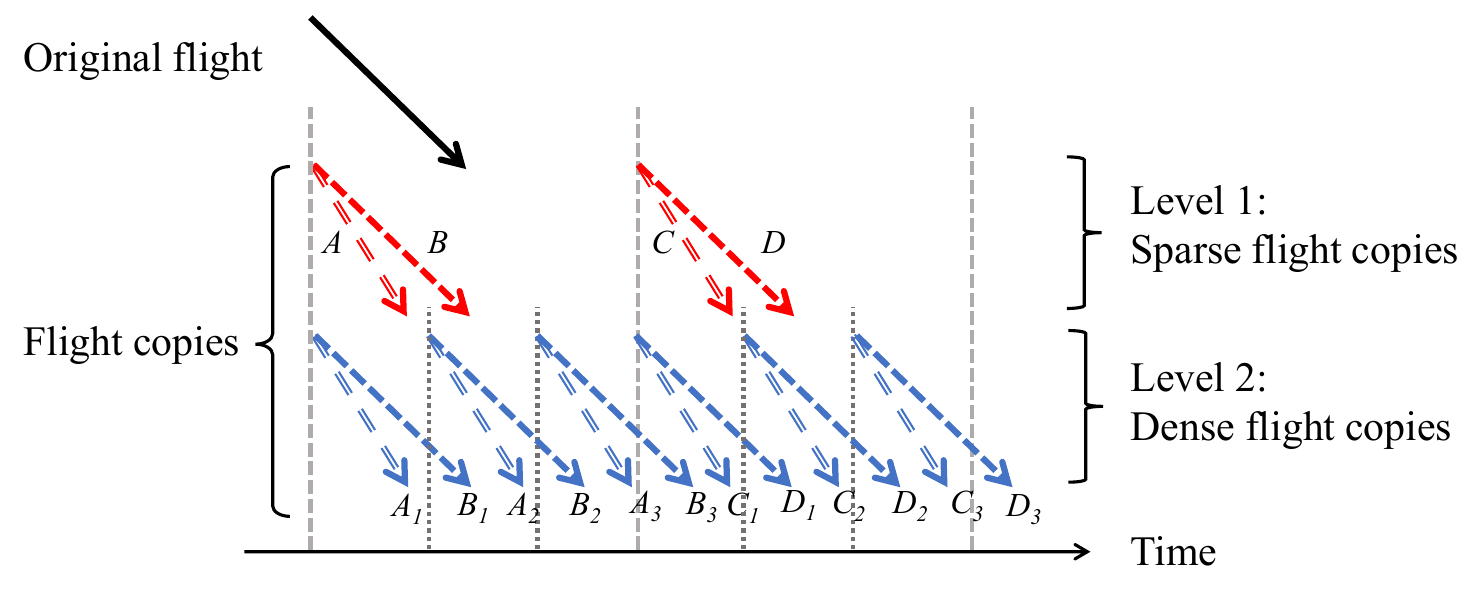}
    \vspace*{-0.5em}
    \caption{Illustration of sparse-dense flight copy approach\label{sparse-denseFlightCopy}}
   \end{minipage}
   \begin{minipage}{1\linewidth}
    \centering
    \vspace*{2em}
    {\includegraphics[width=0.8\textwidth]{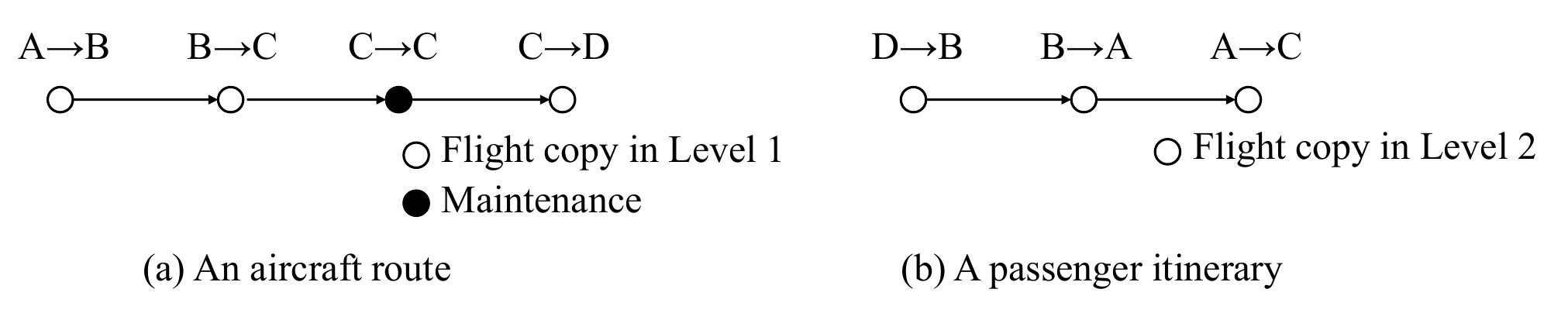}}
    \vspace*{-0.5em}
    \caption{Illustration of an aircraft route and a passenger itinerary\label{routeItinerary}}
   \end{minipage}
\end{figure}

\begin{itemize}
\setlength{\itemsep}{0pt}
\setlength{\parsep}{0pt}
\setlength{\parskip}{0pt}

  \item[(A)] \textit{Recovery time window}: It is a period during which airlines can adjust disrupted schedules, implying that all recovery options must be confined to this period. As is common in the literature on integrated recovery problems, the recovery time windows for different recovery phases are assumed to be the same in this work.

  \item[(B)] \textit{Recovery option}: Various recovery options are available in disruption management. In this paper, the following recovery options are considered: canceling flights, delaying flights, swapping tail assignments, controlling cruise speeds, canceling passenger itineraries, and reallocating passengers. To simplify the integrated recovery problem, swapping tail assignments is restricted to the same fleet type. However, our methodology is also applicable to cases without this restriction.

  \item[(C)] \textit{Slot}: A slot is a period characterized by the following attributes: airport, start time, length of period, and slot capacity. The slot capacity may change when a disruption occurs, which directly affects recovery decisions.

  \item[(D)] \textit{Sparse-dense flight copy}: In airline disruption management, a copy of the original flight is called a \textit{flight copy}. In this paper, we introduce a novel sparse-dense flight copy approach to model flight delays and changes in cruise speeds. As illustrated in Figure \ref{sparse-denseFlightCopy}, based on the length of the copy interval, we divide all flight copies into two levels: sparse flight copies and dense flight copies. The cost of a flight copy consists of the flight delay cost and the cost incurred by cruise speed control (including the cost related to additional fuel consumption and CO$_2$ emission).

  \item[(E)]\textit{Aircraft route}: An aircraft route is a sequence of sparse flight copies performed by an aircraft, as shown in Figure \ref{routeItinerary}(a). The cost of a route is calculated based on the costs incurred by swapping tail assignments and controlling cruise speeds for the flight copies in this route. Specifically, the cost of swapping tail assignments is related to the number of flight copies whose original flights are unplanned for the associated aircraft. The cost of controlling cruise speeds is the sum of the costs incurred by this recovery option for the flight copies contained in the route. The following key constraints are required for aircraft routes.\\
      \hspace*{2em} 1) Match constraint: The route should only contain flight copies performable by the associated aircraft.\\
      \hspace*{2em} 2) Connection constraint: For each two consecutive flight copies in the route, spatial matching and temporal matching (i.e., the minimum turn time requirement) are required. Meanwhile, the start airport and the end airport of the route should be the airports where the aircraft is located at the start time and the end time of the recovery time window, respectively.\\
      \hspace*{2em} 3) Maintenance constraint: The maintenance requirement of an aircraft is regarded as a special sparse flight copy that must be included in the route of that aircraft. Typically, the maintenance station and the required period remain unchanged in airline operations. Hence, a maintenance task is denoted by only one sparse flight copy and one dense flight copy.

  \item[(F)]\textit{Passenger itinerary}: A passenger itinerary is a sequence of dense flight copies for the passengers whose original flight sequences are identical within the recovery time window, as illustrated in Figure \ref{routeItinerary}(b). The cost of a passenger itinerary consists of the passenger delay cost and the reallocation cost. The passenger delay cost is calculated based on the difference between the actual itinerary's end time and the planned itinerary's end time. The reallocation cost is related to the number of unplanned flights in the actual itinerary. Some key constraints for passenger itineraries are as follows.\\
      \hspace*{2em} 1) Connection constraint: A passenger itinerary must also satisfy the connection constraint (spatial-temporal matching and start-end airports), which is the same as that in aircraft routes.\\
      \hspace*{2em} 2) Start time constraint: The start time of a passenger itinerary should not be earlier than the planned start time of the associated original itinerary.
\end{itemize}

\vspace*{0cm}

\subsection{Details of Interactive Mechanism}\label{Section Details of interactive mechanism}
\begin{figure}[!htbp]
\centering
\includegraphics[width=0.9\textwidth]{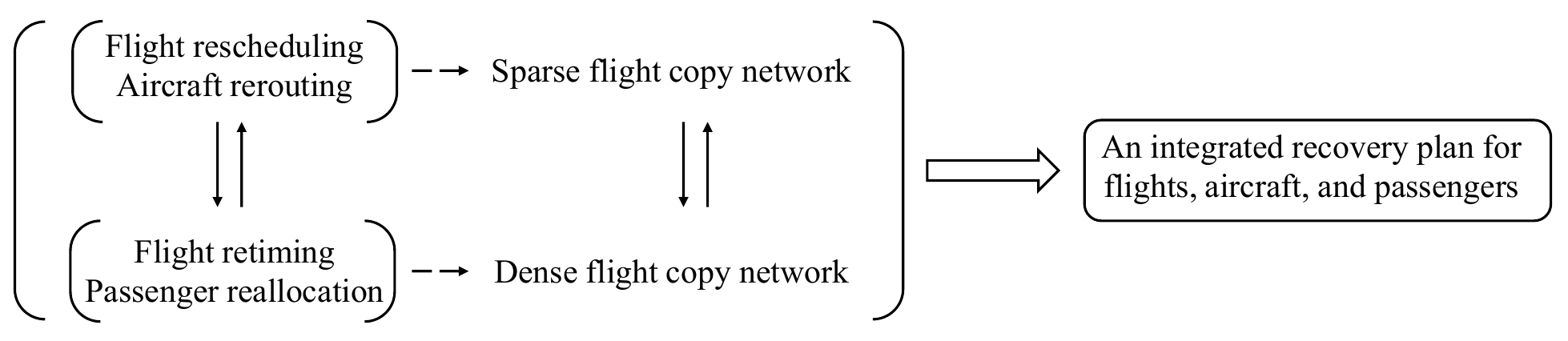}
%\vspace*{-2.2em}
\caption{Illustration of the interactive mechanism\label{interaction}}
\end{figure}

\noindent
After introducing the sparse-dense flight copy approach illustrated in Figure \ref{sparse-denseFlightCopy}, we now develop the sparse-dense flight copy-based interactive mechanism for the integrated flight, aircraft, and passenger recovery problem, which is inspired by the sequential recovery framework in \citeauthor{zhang2016} (\citeyear{zhang2016}). As illustrated in Figure \ref{interaction}, the entire recovery process is divided into two stages. The first stage outputs a rough recovery decision on flight rescheduling and aircraft rerouting, based on the sparse flight copy network generated by flights in Level 1. In the second stage, a slight readjustment of flight departure/arrival times and passenger reallocation are performed on the dense flight copy network generated by flights in Level 2, with passenger reallocation carried out simultaneously. To enable this slight readjustment, we select dense flight copies that are \textit{close} to the sparse flight copies chosen in the first stage. Given a selected sparse flight copy, ``\textit{close}'' means that the flying times of the dense flight copies are the same as that of the sparse flight copy, and the associated delay times fall within a given tolerance relative to the sparse flight copy. For example, dense flight copies $A_1$, $A_2$, and $A_3$ in Figure \ref{sparse-denseFlightCopy} are close to sparse flight copy $A$. Clearly, the problem scales of these two stages are much smaller than that of the integrated recovery problem. However, it is important to note that the recovery decision from the first stage may not be desirable for the integrated recovery, as it does not consider passenger reallocation when making the decision on flight rescheduling and aircraft rerouting. To overcome this issue, the decision from the first stage should be updated based on feedback from the second stage. Hence, the final integrated recovery plan is obtained by alternately solving these two stages, addressing flight recovery, aircraft rerouting, and passenger reallocation from an integrated perspective.

We use the flight copies in Figure \ref{sparse-denseFlightCopy} as an example to demonstrate the process of the proposed interactive mechanism. Recall that aircraft routes and passenger itineraries are generated using sparse flight copies and dense flight copies, respectively (Section \ref{Section Definitions} (E)$-$(F)). In the first stage, we obtain a decision on flight and aircraft recovery by selecting a sparse flight copy for the original flight and assigning it to an aircraft. Suppose we select sparse flight copy $A$ in this stage. In the second stage, we readjust the flight recovery decision and reallocate passengers by choosing a dense flight copy from $A_1$, $A_2$, and $A_3$, rather than from all dense flight copies. If the performance of the resulting integrated recovery plan is not satisfactory, we feed this information back and update the previous decision in the first stage by selecting a new sparse flight copy from $B$, $C$, or $D$. The interaction between the two stages is repeated until a desirable integrated recovery plan is obtained.

The advantages of the proposed interactive mechanism are twofold. On the one hand, rather than solving the integrated recovery problem sequentially, the interactive mechanism captures the interdependencies between different recovery phases by alternately adjusting their decisions and exchanging associated feedback. This represents one of the main improvements over the sequential recovery framework in \citeauthor{zhang2016} (\citeyear{zhang2016}). On the other hand, with the help of the sparse-dense flight copy approach, the interactive mechanism reduces solution difficulty by dividing the integrated recovery problem into two stages whose problem scales are far smaller than that of the conventional flight copy approach. We again use Figure \ref{sparse-denseFlightCopy} as an example. By using sparse flight copies, the workload of obtaining a flight and aircraft recovery plan is greatly reduced. Benefiting from selecting dense flight copies close to the chosen sparse flight copy, we only need to consider three dense flight copies in the second stage instead of all twelve. This advantage becomes even more evident in practical airline operations, where a large number of flight copies are needed to represent flight delays and cruise speeds on a large-scale airline network.

\begin{figure}[!htbp]
\centering
\hspace*{-0.4em}
\includegraphics[width=1\textwidth]{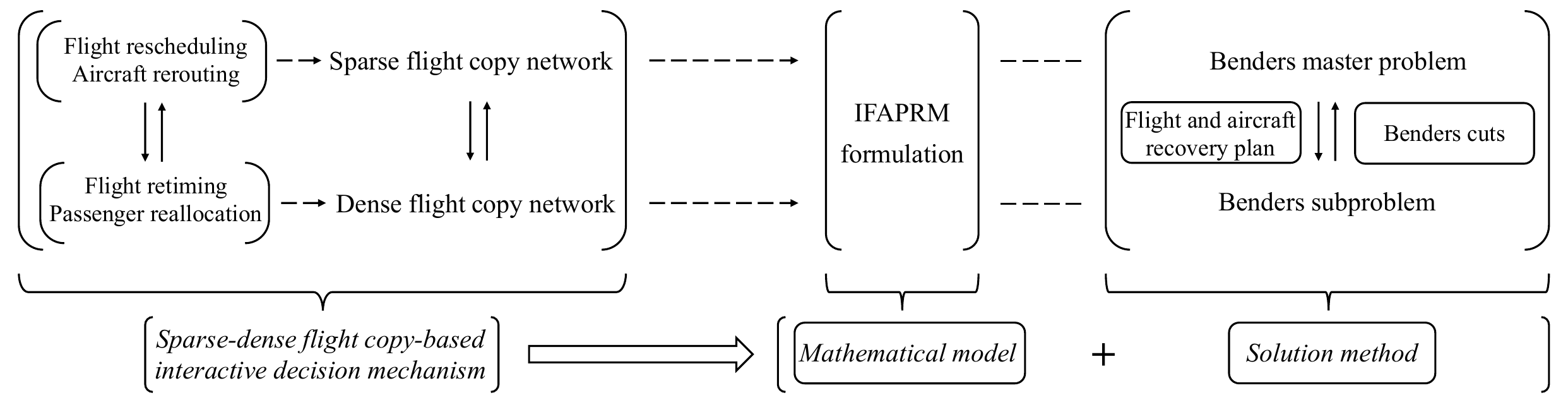}
\caption{Operationalization of the interactive mechanism\label{interaction2}}
\end{figure}

A central task, then, is how to make the interactive mechanism operational, particularly how to update decisions in the sparse flight copy network based on interactive feedback from the dense flight copy network. Instead of directly adopting an alternating minimization method, which struggles to effectively coordinate the interdependent decisions between the two networks, we develop an integrated recovery model and propose a customized Benders decomposition method to operationalize this mechanism. As illustrated in Figure \ref{interaction2}, the developed model corresponds to the ``{\it Mathematical model}'' component, while the proposed method corresponds to the ``{\it Solution method}'' component. These are elaborated in Sections \ref{Section Mathematical model} and \ref{Section Solution methodology}, respectively.

\section{Mathematical Model}\label{Section Mathematical model}
\noindent
As a critical component of operationalizing the interactive mechanism, an integrated flight, aircraft, and passenger recovery model (IFAPRM) is developed in this section to formulate the integrated recovery problem.\par

\subsection{Notations}\label{Section Notations}
\noindent
Before providing the specific formulation, we give the notations used in this paper, which is summarized in Table \ref{Notations}.\par

\begin{center}
\small
\renewcommand{\baselinestretch}{1}
    \vspace*{-0.5cm}
	\begin{longtable}{p{0.15\textwidth}p{0.78\textwidth}}
        \caption{Notations used in this paper}\label{Notations}
		\endfirsthead
			
        \midrule
		\endhead
			
        \midrule
        \multicolumn{2}{l}{\footnotesize Continued on next page}
		\endfoot
			
		\hline
	\endlastfoot

    \toprule
    \multicolumn{2}{l}{\bf{Sets}}\\
    \midrule
    $F$& the set of original flights\\
    $FC_f^S$& the set of the sparse flight copies for flight~$f\in F$\\
    $FC_{f}^D$& the set of the dense flight copies for flight~$f\in F$\\
    $FC_{f_u}^D$& the set of the dense flight copies for sparse flight~$f_u\in FC_f^S$\\
    $R$& the set of aircraft\\
    $P^r$& the set of the feasible routes for aircraft $r\in R$\\
    $\bar{P}^r$ &the subset of the feasible routes for aircraft $r\in R$, i.e., $\bar{P}^r \subseteqq P^r$ \\
    $E^R$& the set of the flight connections that can be performed by an aircraft\\
    $OI$& the set of original itineraries\\
    $P^i$& the set of the feasible itineraries for original itinerary $i\in OI$\\
    $\bar{P}^i$ &the subset of the feasible itineraries for original itinerary $i\in OI$, i.e., $\bar{P}^i \subseteqq P^i$ \\
    $SL$& the set of slots\\
    \toprule
    \multicolumn{2}{l}{\bf{Parameters}}\\
    \midrule
    $d_{f}^{F}$ & the cancellation cost of flight $f$\\
    $c_{f}^{OF}$ & the original fuel cost of flight $f$\\
    $d_{i}^{OI}$& the unit itinerary cancellation cost for each passenger from original itinerary $i \in OI$\\
    $c_{i}^{DT}$& the unit delay cost for each passenger from original itinerary $i \in OI$\\
    $c_{i}^{C}$& the unit change cost for each passenger from original itinerary $i \in OI$ if the associated flight sequence changes\\
    $c_{f_{uv}}$& the cost of dense flight copy $f_{uv}$\\
    $c_p^r$& the cost of aircraft route $p\in P^r$\\
    $n_{i}^{OI}$& the number of the passengers in original itinerary $i \in OI$\\
    $Cap_r$& the number of the seats in aircraft $r$\\
    $MaxCap$& the maximum number of seats among all aircraft\\
    $MinTR^{r}_{(f_1,f_2)}$& the minimum turn time for connection $(f_1,f_2)$ if it is performed by aircraft $r \in R$\\
    $MaxT$& the end time of the recovery time window\\
    $T_{f_{uv}}^{dep}/T_{f_{uv}}^{arr}$& the departure/arrival time of dense flight copy $f_{uv}$\\
    $DT_{p}^{i}$& the delay time of itinerary $p\in P^i$ related to its original itinerary $i \in OI$\\\
    $u_{sdep}/u_{sarr}$& the departure/arrival capacity of slot $s\in SL$\\
    $a_{{p,{f_u}}}^{r}$ &1, if sparse flight copy $f_u$ is in route $p\in P^r$;~0, otherwise\\
    $b_{p,{(f_1,f_2)}}^{r}$ &1, if flight connection $(f_1,f_2)$ is in route $p\in P^r$;~0, otherwise\\
    $a_{{p,{f_{uv}}}}^{i}$ &1, if dense flight copy $f_{uv}$ is in itinerary $p\in P^i$;~0, otherwise\\
    $a_{f_{uv}}^{sdep}/a_{f_{uv}}^{sarr}$ &1, if dense flight copy $f_{uv}$ departs/arrives within slot $s \in SL$;~0, otherwise\\
    $n_{i,p}^{C}$& the number of unplanned flights in itinerary $p\in P^i$ for original itinerary $i \in OI$\\
    \toprule
    \multicolumn{2}{l}{\bf{Decision Variables}}\\
    \midrule
    $z_{f}^{F}$ &1, if original flight~$f\in F$ is canceled;~0, otherwise\\
    $y_p^r$ &1, if aircraft $r\in R$ is assigned to route $p\in P^r$;~0, otherwise\\
    $x_{f_{uv}}$ &1, if dense flight copy $f_{uv}$ is chosen by flight $f \in F$;~0, otherwise\\
    $z_i^{OI}$& the number of the passengers who are from original itinerary $i\in OI$ and unassigned\\
    $w_p^{i}$& the number of the passengers who are from original itinerary $i\in OI$ and assigned to itinerary $p\in P^i$\\
    \toprule
    \multicolumn{2}{l}{\bf{Main Abbreviations}}\\
    \midrule
    IFAPRM & the integrated flight, aircraft, and passenger recovery model $\eqref{IFAPRMobj}-\eqref{IFAPRMcon14}$\\
    FRARM & the flight rescheduling with aircraft rerouting model $\eqref{ARMobj}-\eqref{ARMcon7}$\\
    FRPRM & the flight retiming with passenger reallocation model $\eqref{FRPRMobj}-\eqref{FRPRMcon10}$\\
    BD & Benders decomposition\\
    BMP & Benders master problem\\
    RxBMP & the relaxed BMP, i.e., the BMP with the partial Benders cuts\\
    LR-RxBMP & the linear relaxation of RxBMP\\
    RtRxBMP & the restricted RxBMP, i.e., the RxBMP with partial routes $\bar{P}^{r}$\\
    LR-RtRxBMP & the linear relaxation of RtRxBMP\\
    BSP & Benders subproblem\\
    LR-BSP & the linear relaxation of BSP\\
    RtBSP & the restricted BSP, i.e., the BSP with partial itineraries $\bar{P}^{i}$\\
    LR-RtBSP & the linear relaxation of RtBSP\\
    P\uppercase\expandafter{\romannumeral1}LR-RtBSP & the phase \uppercase\expandafter{\romannumeral1} optimization problem of LR-RtBSP\\
    \bottomrule
	\end{longtable}
    \vspace*{-0.5cm}
\end{center}

\subsection{Formulation}\label{Section Formulation}
\noindent
Various formulations have been developed for the integrated recovery problem in previous studies (\citeauthor{sinclair2014}, \citeyear{sinclair2014}; \citeauthor{zhang2016}, \citeyear{zhang2016}; \citeauthor{marla2017}, \citeyear{marla2017}; etc). Some studies used an arc-based modelling method to formulate the problem. While this approach can explicitly describe the practical restrictions, it requires a large number of variables and constraints, which imposes a significant computational burden on solving the integrated recovery problem. To overcome this issue, we apply a path-based modelling method and use the column generation technique to alleviate the workload of solving the problem.\par

The IFAPRM formulation consists of four parts: flight rescheduling, aircraft rerouting, flight retiming, and passenger reallocation. The formulation for flight rescheduling and aircraft rerouting is a classical integrated flight and aircraft recovery model, as seen in \citeauthor{Rosenberger2003} (\citeyear{Rosenberger2003}) and \citeauthor{liang2018} (\citeyear{liang2018}). To operationalize the proposed interactive mechanism, the latter two parts are formulated based on the decision of the former two. Specifically, the integrated recovery problem is formulated as follows.\par

\vspace*{-0.8cm}
\begin{eqnarray}
\nonumber
&&\hspace{-1.2cm}\textbf{IFAPRM formulation:}\\
\nonumber
&\min &\sum_{f\in F} (d_{f}^{F} - c_{f}^{OF})z_{f}^{F} + \sum_{r\in R}\sum_{p\in {P^r}}c_p^ry_p^r +
\sum_{f \in F} \sum_{f_{u} \in {FC_{f}^S}} \sum_{f_{uv} \in {FC_{f_{u}}^D}} c_{f_{uv}} x_{f_{uv}} \label{IFAPRMobj}\\
&&
 + \sum_{i \in OI}d_{i}^{OI}z_{i}^{OI}+\sum_{i \in OI}\sum_{p \in P^i}c_{i}^{DT}DT_{p}^{i}w_{p}^{i} + \sum_{i \in OI}\sum_{p \in P^i}c_{i}^{C}n_{i,p}^{C}w_{p}^{i}\\
&\mathrm{s.t.} &\sum_{r\in R}\sum_{p\in {P^r}}\sum_{f_u\in {FC_f^S}}a_{{p,{f_u}}}^{r} y_p^r+z_{f}^{F}=1,\ \ \forall f\in {F},\label{IFAPRMcon1}\\
&&
\sum_{p\in {P^r}} y_p^r \leq 1,\ \ \forall r\in R,\label{IFAPRMcon2}\\
&&
\sum_{f_{uv} \in {FC_{f_u}^D}} x_{f_{uv}} = \sum_{r\in R}\sum_{p\in {P^r}} a_{{p,{f_u}}}^{r} y_p^r,\ \ \forall f\in {F}, f_u \in {FC_{f}^S},\label{IFAPRMcon3}\\
\nonumber
&&
\sum_{f_{1u} \in {FC_{f_{1}}^S}} \sum_{f_{1uv} \in {FC_{f_{1u}}^D}} T_{f_{1uv}}^{arr} x_{f_{1uv}} + \sum_{r \in R} \sum_{p\in {P^r}} (MinTR^{r}_{(f_1,f_2)}+MaxT)b_{p,{(f_1,f_2)}}^{r}  y_p^r\\
&&\hspace{1cm} \leq \sum_{f_{2u} \in {FC_{f_{2}}^S}} \sum_{f_{2uv} \in {FC_{f_{2u}}^D}} T_{f_{2uv}}^{dep} x_{f_{2uv}} + MaxT,\ \ \forall (f_1, f_2) \in {E^R},\label{IFAPRMcon4}\\
&&
\sum_{f \in F} \sum_{f_{u} \in {FC_{f}^S}} \sum_{f_{uv} \in {FC_{f_{u}}^D}} a_{f_{uv}}^{sdep} x_{f_{uv}} \leq u_{sdep},\ \ \forall s\in {SL},\label{IFAPRMcon5}\\
&&
\sum_{f \in F} \sum_{f_{u} \in {FC_{f}^S}} \sum_{f_{uv} \in {FC_{f_{u}}^D}} a_{f_{uv}}^{sarr} x_{f_{uv}} \leq u_{sarr},\ \ \forall s\in {SL},\label{IFAPRMcon6}\\
&&
\sum_{p\in {P^i}} w_p^i + z_i^{OI} = n_{i}^{OI},\ \ \forall i\in OI,\label{IFAPRMcon7}\\
&&
\sum_{i \in OI}\sum_{p \in P^i} a_{{p,{f_{uv}}}}^{i} w_p^i \leq MaxCap \cdot x_{f_{uv}},\ \ \forall f \in F, f_{u} \in {FC_{f}^S}, f_{uv} \in {FC_{f_{u}}^D},\label{IFAPRMcon8}\\
\nonumber
&&
\sum_{i \in OI}\sum_{p \in P^i}\sum_{f_{u} \in {FC_{f}^S}} \sum_{f_{uv} \in {FC_{f_u}^D}} a_{{p,{f_{uv}}}}^{i}w_{p}^{i}
\leq \\
&&
\hspace{3cm} \sum_{r \in R}Cap_r \sum_{p \in P^r}\sum_{f_{u} \in FC_{f}^S}a_{{p,{f_u}}}^{r} y_p^r + MaxCap\cdot z_{f}^{F},\ \ \forall f \in F,\label{IFAPRMcon9}\\
&&
z_{f}^{F}\in \{0,1\},\ \ \forall j\in {F},\label{IFAPRMcon10}\\
&&
y_p^r\in \{0,1\},\ \ \forall r\in R, p\in P^r,\label{IFAPRMcon11}\\
&&
x_{f_{uv}}\in \{0,1\},\ \ \forall f \in F, f_{u} \in FC_{f}^S, f_{uv} \in FC_{f_{u}}^D,\label{IFAPRMcon12}\\
&&
z_i^{OI} \in \mathbb{Z_{+}},\ \ \forall i\in OI,\label{IFAPRMcon13}\\
&&
w_p^{i} \in \mathbb{Z_{+}},\ \ \forall i\in OI, p \in P^i.\label{IFAPRMcon14}
\end{eqnarray}

The objective function \eqref{IFAPRMobj} is to minimize the total recovery cost, which includes the costs of canceling flights, assigning aircraft routes, delaying flights, canceling passenger itineraries, delaying passenger arrivals, and changing passenger itineraries. The calculation of cost parameters $c_{f_{uv}},~c_{p}^{r}$ and $c_{i}^{DT}$ has been discussed in Section \ref{Section Definitions} (D)$-$(F). In this work, the cost incurred by cruise speed control is calculated by the functions in \citeauthor{akturk2014} (\citeyear{akturk2014}). Constraints \eqref{IFAPRMcon1} ensure that each flight is either operated by an aircraft or canceled. Constraints \eqref{IFAPRMcon2} mean that each aircraft is assigned at most one route. Constraints \eqref{IFAPRMcon3} guarantee that each flight's recovery option (delay or cancellation) must be consistent between aircraft rerouting and flight retiming. Constraints \eqref{IFAPRMcon4} ensure that the minimum turn time requirement is satisfied for each flight connection in the chosen aircraft routes after retiming flights. Constraints \eqref{IFAPRMcon5} and \eqref{IFAPRMcon6} capture the departure and arrival capacity restrictions for slots, respectively. Constraints \eqref{IFAPRMcon7} state that each passenger is either assigned to an itinerary or left unassigned during disruptions. Constraints \eqref{IFAPRMcon8} ensure that passengers can only be reallocated to the chosen dense flight copies. Constraints \eqref{IFAPRMcon9} are the seat capacity restrictions for each flight.\par

The following observations can be drawn from the above model. On the one hand, although the IFAPRM formulation is complex due to the presence of the two types of flight copies, the sparse-dense flight copy approach naturally introduces a two-stage structure. This structure enables the application of Benders decomposition, which divides the integrated recovery model into a Benders master problem on the sparse flight copy network and a Benders subproblem on the dense flight copy network (see Figure \ref{interaction2}). Compared with the model formulated by the conventional flight copy approach, the problem scales of these two problems are far smaller, which greatly reduces the solution difficulty. On the other hand, with the path-based modelling method, column generation can be employed to further reduce solution difficulty by avoiding the enumeration of all possible aircraft routes and passenger itineraries.

\section{Solution Method}\label{Section Solution methodology}
\noindent
As illustrated in Figure \ref{interaction2}, the operationalization of the proposed interactive mechanism consists of the IFAPRM formulation and the corresponding solution method. For the integrated recovery problem on a practical airline network, IFAPRM is a large-scale integer programming problem with complex constraints, so a satisfactory recovery plan cannot be directly obtained from a commercial solver within reasonable runtime. By exploiting the inherent multi-stage characteristic of the IFAPRM formulation, we adopt a Benders decomposition strategy to solve the integrated model. Notably, the iteration between Benders master problem and Benders subproblem in the Benders decomposition framework (\citeauthor{benders1962}, \citeyear{benders1962}) coincides perfectly with the iterative adjustment in the interactive mechanism. This alignment establishes a bidirectional feedback loop that enables the feedback from one stage to guide decision-making in the other, leading to higher-quality recovery solutions. In addition, to avoid enumerating all possible aircraft routes and passenger itineraries, column generation (\citeauthor{dantzig1960}, \citeyear{dantzig1960}) is applied to generate routes and itineraries iteratively. In this section, we first present the Benders subproblem and the Benders master problem for the IFAPRM formulation. Then, we illustrate the specific process of the customized Benders decomposition (CBD) method.

\subsection{Benders Subproblem}\label{Section Benders subproblem}
\noindent
According to the interactive mechanism, flight retiming and passenger reallocation are performed after a decision on flight and aircraft recovery is made. Thus, the problem that contains flight retiming and passenger reallocation is taken as Benders subproblem (BSP) in this paper, which consists of the last four terms in the objective function \eqref{IFAPRMobj}, Constraints $\eqref{IFAPRMcon3}-\eqref{IFAPRMcon9}$, and Constraints $\eqref{IFAPRMcon12}-\eqref{IFAPRMcon14}$.\par

Given a solution $(\bar{y}_p^r, \bar{z}_{f}^{F})$ corresponding to a decision of the flight scheduling with aircraft rerouting, the associated BSP is presented as follows, which is a flight retiming with passenger reallocation model (FRPRM) for the solution $(\bar{y}_p^r, \bar{z}_{f}^{F})$.\par

\vspace{-1em}
\begin{eqnarray}
&&\hspace{-1.5cm}\textbf{FRPRM formulation:}\notag\\
&\min &\sum_{f \in F} \sum_{f_{u} \in {FC_{f}^S}} \sum_{f_{uv} \in {FC_{f_{u}}^D}} c_{f_{uv}} x_{f_{uv}} + \sum_{i \in OI}d_{i}^{OI}z_{i}^{OI}+\sum_{i \in OI}\sum_{p \in P^i}c_{i}^{DT}DT_{p}^{i}w_{p}^{i}\notag\\
&&\hspace{0cm} + \sum_{i \in OI}\sum_{p \in P^i}c_{i}^{C}n_{i,p}^{C}w_{p}^{i}\label{FRPRMobj}\\
&\mathrm{s.t.} &\sum_{f_{uv} \in {FC_{f_u}^D}} x_{f_{uv}} = \sum_{r\in R}\sum_{p\in {P^r}} a_{{p,{f_u}}}^{r} \bar{y}_p^r,\ \ \forall f\in {F}, f_u \in {FC_{f}^S},\label{FRPRMcon1}\\
\nonumber
&&
\sum_{f_{1u} \in {FC_{f_{1}}^S}} \sum_{f_{1uv} \in {FC_{f_{1u}}^D}} T_{f_{1uv}}^{arr} x_{f_{1uv}} + \sum_{r \in R} \sum_{p\in {P^r}} (MinTR^{r}_{(f_1,f_2)}+MaxT)b_{p,{(f_1,f_2)}}^{r}  \bar{y}_p^r\\
&&\hspace{1cm} \leq \sum_{f_{2u} \in {FC_{f_{2}}^S}} \sum_{f_{2uv} \in {FC_{f_{2u}}^D}} T_{f_{2uv}}^{dep} x_{f_{2uv}} + MaxT,\ \ \forall (f_1, f_2) \in {E^R},\label{FRPRMcon2}\\
&&
\sum_{f \in F} \sum_{f_{u} \in {FC_{f}^S}} \sum_{f_{uv} \in {FC_{f_{u}}^D}} a_{f_{uv}}^{sdep} x_{f_{uv}} \leq u_{sdep},\ \ \forall s\in {SL},\label{FRPRMcon3}\\
&&
\sum_{f \in F} \sum_{f_{u} \in {FC_{f}^S}} \sum_{f_{uv} \in {FC_{f_{u}}^D}} a_{f_{uv}}^{sarr} x_{f_{uv}} \leq u_{sarr},\ \ \forall s\in {SL},\label{FRPRMcon4}\\
&&
\sum_{p\in {P^i}} w_p^i + z_i^{OI} = n_{i}^{OI},\ \ \forall i\in OI,\label{FRPRMcon5}\\
&&
\sum_{i \in OI}\sum_{p \in P^i} a_{{p,{f_{uv}}}}^{i} w_p^i \leq MaxCap \cdot x_{f_{uv}},\ \ \forall f \in F, f_{u} \in {FC_{f}^S}, f_{uv} \in {FC_{f_{u}}^D},\label{FRPRMcon6}\\
\nonumber
&&
\sum_{i \in OI}\sum_{p \in P^i}\sum_{f_{u} \in {FC_{f}^S}} \sum_{f_{uv} \in {FC_{f_u}^D}} a_{{p,{f_{uv}}}}^{i}w_{p}^{i}
\leq \\
&&
\hspace{3cm} \sum_{r \in R}Cap_r \sum_{p \in P^r}\sum_{f_{u} \in FC_{f}^S}a_{{p,{f_u}}}^{r} \bar{y}_p^r + MaxCap\cdot \bar{z}_{f}^{F},\ \ \forall f \in F,\label{FRPRMcon7}\\
&&
x_{f_{uv}}\in \{0,1\},\ \ \forall f \in F, f_{u} \in FC_{f}^S, f_{uv} \in FC_{f_{u}}^D,\label{FRPRMcon8}\\
&&
z_i^{OI} \in \mathbb{Z_{+}},\ \ \forall i\in OI,\label{FRPRMcon9}\\
&&
w_p^{i} \in \mathbb{Z_{+}},\ \ \forall i\in OI, p \in P^i.\label{FRPRMcon10}
\end{eqnarray}

The Benders cuts of the BSP are obtained from the linear relaxation of the BSP (LR-BSP). Let $\Pi_{Fea}$ and $\Pi_{Opt}$ denote the sets of extreme directions and extreme points in the dual feasible region of the LR-BSP, respectively, and let
\begin{eqnarray*}
\pi := (\pi_{f_u}^{1}, \pi_{(f_1, f_2)}^{2}, \pi_{s}^{3}, \pi_{s}^{4}, \pi_{i}^{5}, \pi_{f_{uv}}^{6}, \pi_{f}^{7}) \label{the dual values of LR-BSP}
\end{eqnarray*}
be the vector of dual values related to Constraints $\eqref{FRPRMcon1}-\eqref{FRPRMcon7}$ in the LR-BSP. Then, the Benders feasibility cuts and Benders optimality cuts are
\begin{eqnarray}
\sum_{r\in R}\sum_{p\in {P^r}}a_{p}^{r}(\pi)y_p^r  + \sum_{f \in F}b_{f}(\pi)z_{f}^{F} + c(\pi) \leq 0,
\ \ \forall \pi \in \Pi_{Fea}, \label{Benders Fea Cut}\\
\sum_{r\in R}\sum_{p\in {P^r}}a_{p}^{r}(\pi)y_p^r  + \sum_{f \in F}b_{f}(\pi)z_{f}^{F} + c(\pi)\leq q,
\ \ \forall \pi \in \Pi_{Opt}, \label{Benders Opt Cut}
\end{eqnarray}
where $a_{p}^{r}(\pi),~b_{f}(\pi)$ and $c(\pi)$ are defined as
\begin{small}
\begin{eqnarray*}
&a_{p}^{r}(\pi) &= \sum\limits_{f \in F} \sum\limits_{f_{u} \in {FC_{f}^S}}\pi_{f_u}^{1} a_{{p,{f_u}}}^{r} + \sum\limits_{(f_1, f_2) \in {E^R}} \pi_{(f_1, f_2)}^{2}(MinTR^{r}_{(f_1,f_2)}+MaxT)b_{p,{(f_1,f_2)}}^{r}\\
&&\hspace{0.5cm}  - \sum\limits_{f \in F} \sum\limits_{f_{u} \in {FC_{f}^S}}\pi_{f}^{7}Cap_r a_{{p,{f_u}}}^{r},\\
&b_{f}(\pi) &= - \pi_{f}^{7}MaxCap,\\
&c(\pi) &= -\sum\limits_{(f_1, f_2) \in {E^R}}\pi_{(f_1, f_2)}^{2}MaxT -\sum\limits_{s \in SL}(\pi_{s}^{3}u_{sdep}+\pi_{s}^{4}u_{sarr})
 + \sum\limits_{i \in OI} \pi_{i}^{5}n_{i}^{OI},
\end{eqnarray*}
\end{small}\par
\noindent and $q$ is a decision variable linking the Benders master problem and the Benders subproblem in Benders decomposition framework.\par

Considering the large scale of the passenger itineraries in $P^{i}$, we apply column generation technique to generate itineraries gradually. New elements in $P^{i}$ are generated by iteratively solving the CG restricted master problem and the CG subproblem. The CG restricted master problem of the BSP is the linear relaxation of the restricted BSP (whose passenger itineraries are restricted to the set $\bar{P}^{i} \subseteq P^{i}$). Given the dual solutions of the CG restricted master problem, the CG subproblem aims to find passenger itineraries with the most negative reduced cost for each original itinerary $i \in OI$ on the associated dense flight copy connection network (similar to the aircraft connection network in \citeauthor{papadakos2009} (\citeyear{papadakos2009})), where the reduced cost $\bar{c}_p^{i}$ of itinerary $p \in P^{i}$ is
\begin{eqnarray}
\bar{c}_p^{i} := c_{i}^{DT} DT_{p}^{i} + c_{i}^{C}n_{i,p}^{C} - \pi_{i}^{5}+ \sum_{f \in F} \sum_{f_{u} \in {FC_{f}^S}} \sum_{f_{uv} \in {FC_{f_{u}}^D}}(\pi_{f_{uv}}^{6}+\pi_{f}^{7})a_{{p,{f_{uv}}}}^{i}. \label{reduced cost itinerary}
\end{eqnarray}
It is noted that $DT_{p}^{i}, n_{i,p}^{C}$ and $a_{{p,{f_{uv}}}}^{i}$ are decision variables in the CG subproblem. In this paper, the CG subproblem is solved by the constraint programming (CP) solver of CPLEX, whose detailed model is given in Appendix B.1.\par

Since the BSP is solved by the CG process, the feasibility and optimality of the BSP are not equivalent to those of the RtBSP. The specific processes of checking the feasibility and the optimality are summarized by the \textit{\textbf{BSP Fea}} block and the \textit{\textbf{BSP Opt}} block in Figure \ref{CBDflow}, respectively. The \textit{phase \uppercase\expandafter{\romannumeral1} optimization problem} (\citeauthor{boyd2004}, \citeyear{boyd2004}) of the linear relaxation of the RtBSP (P\uppercase\expandafter{\romannumeral1}LR-RtBSP) is introduced to help us judge the feasibility of the BSP when the RtBSP is infeasible. In addition, by exploiting the structure of the BSP, we further develop several acceleration techniques for the CBD method, which will be discussed in Section \ref{Section Acceleration techniques}.\par

\subsection{Benders Master Problem}\label{Section Benders master problem}
\noindent
The Benders master problem (BMP) is a flight rescheduling with aircraft rerouting model (FRARM), which provides a decision for the flight and aircraft recovery. For FRARM, the objective function consists of the first two terms of \eqref{IFAPRMobj} plus a decision variable $q$. The constraints of FRARM includes Constraints $\eqref{IFAPRMcon1}-\eqref{IFAPRMcon2}$, $\eqref{IFAPRMcon10}-\eqref{IFAPRMcon11}$, the Benders cuts $\eqref{Benders Fea Cut}-\eqref{Benders Opt Cut}$, and a nonnegative constraint for $q$. Hence, the BMP is given as follows.\par

\vspace{-1em}
\begin{eqnarray}
&&\hspace{-2cm}\textbf{FRARM formulation:}\notag\\
&\min &\sum_{f\in F}(d_{f}^{F} - c_{f}^{OF})z_{f}^{F} + \sum_{r\in R}\sum_{p\in {P^r}}c_p^ry_p^r + q \label{ARMobj}\\
&\mathrm{s.t.} &\sum_{r\in R}\sum_{p\in {P^r}}\sum_{f_u\in {FC_f^S}}a_{{p,{f_u}}}^{r} y_p^r+z_{f}^{F}=1,\ \ \forall f\in {F},\label{ARMcon1}\\
&&
\sum_{p\in {P^r}} y_p^r \leq 1,\ \ \forall r\in R,\label{ARMcon2}\\
&&
\sum_{r\in R}\sum_{p\in {P^r}}a_{p}^{r}(\pi)y_p^r  + \sum_{f \in F}b_{f}(\pi)z_{f}^{F} + c(\pi) \leq 0,
\ \ \forall \pi \in \Pi_{Fea}, \label{ARMcon3}\\
&&
\sum_{r\in R}\sum_{p\in {P^r}}a_{p}^{r}(\pi)y_p^r  + \sum_{f \in F}b_{f}(\pi)z_{f}^{F} + c(\pi)\leq q,
\ \ \forall \pi \in \Pi_{Opt}, \label{ARMcon4}\\
&&
z_{f}^{F}\in \{0,1\},\ \ \forall j\in {F},\label{ARMcon5}\\
&&
y_p^r\in \{0,1\},\ \ \forall r\in R, p\in P^r,\label{ARMcon6}\\
&&
q\geq 0.\label{ARMcon7}
\end{eqnarray}\par

It is known that the sizes of $\Pi_{Fea}$ and $\Pi_{Opt}$ grow exponentially with the size of input parameters, so enumerating all their elements is impractical. Besides, most constraints corresponding to the elements in $\Pi_{Fea}$ and $\Pi_{Opt}$ are not active at the optimal solution. To overcome the above two issues, we apply an iterative way to gradually find useful elements in these sets. In the iteration process of the CBD method, Constraints \eqref{ARMcon3} and \eqref{ARMcon4} are relaxed to the constraints related to the sets $\overline{\Pi}_{Fea} \subseteq \Pi_{Fea}$ and $\overline{\Pi}_{Opt} \subseteq \Pi_{Opt}$ respectively, i.e., the BMP turns to be a relaxed Benders master problem (RxBMP). The elements in $\Pi_{Fea}$ and $\Pi_{Opt}$ are iteratively generated by solving the BSP and then added to $\overline{\Pi}_{Fea}$ and $\overline{\Pi}_{Opt}$.\par

Similar to the passenger itinerary set $P^{i}$ in the previous subsection, it is unrealistic to enumerate all aircraft routes in $P^r$. Hence, we also apply the CG technique to solve the RxBMP. During the CG process, the aircraft routes set $P^{r}$ is restricted to the set $\bar{P}^{r} \subseteq P^{r}$, i.e., the RxBMP becomes the restricted RxBMP (RtRxBMP). The CG restricted master problem of the RxBMP is the linear relaxation of the RtRxBMP. The associated CG subproblem is a resource-constrained shortest path problem on a sparse flight copy connection network, which aims to find the aircraft routes with the most negative reduced cost. For aircraft $r \in R$, the reduced cost $\bar{c}_p^{r}$ of aircraft route $p \in P^r$ is
\begin{eqnarray}
 \bar{c}_p^{r} := c_p^{r}- \sum\limits_{f\in {F}}\sum\limits_{f_u\in {FC_f^S}}\lambda_{f}^{1}a_{{p,{f_u}}}^{r} + \lambda_{r}^{2} + \sum\limits_{\pi\in {\overline{\Pi}_{Fea}}} \lambda_{\pi}^{3}a_{p}^{r}(\pi) +  \sum\limits_{\pi\in {\overline{\Pi}_{Opt}}} \lambda_{\pi}^{4}a_{p}^{r}(\pi), \label{reduced cost route}
\end{eqnarray}
where $\lambda_{f}^{1}, \lambda_{r}^{2}, \lambda_{\pi}^{3}$ and $\lambda_{\pi}^{4}$ are the dual values of Constraints $\eqref{ARMcon1}-\eqref{ARMcon4}$. The associated CG subproblem is also solved by the CP solver of CPLEX in this paper, whose detailed model is given in Appendix B.2.

The RxBMP is always feasible because the vector
\begin{eqnarray*}
({y_{p}^{r}}, {z_{f}^{F}}, q) := \left({\bf{0}}, {\bf{1}}, \max\limits_{\pi \in \overline{\Pi}_{Opt}} \left\{\sum\limits_{f \in F} b_{f}(\pi) + c(\pi)\right\}\right)
\end{eqnarray*}
is a feasible solution for the RtRxBMP (the restricted problem of RxBMP). Hence, we only need to use CG to solve the RxBMP for its optimality, which is illustrated by the \textit{\textbf{RxBMP Opt}} block in Figure \ref{CBDflow}. In this figure, $v_{LM}(\bar{y}_p^r, \bar{z}_f^F, \bar{q})$ denotes the optimal value of the linear relaxation of the current RtRxBMP.

The BMP aims to output a decision on the flight and aircraft recovery with the help of the Benders cuts. Furthermore, we will introduce some valid inequalities in Section \ref{Section Valid inequalities} to tighten the formulation of the BMP.\par

\subsection{A Customized Benders Decomposition Method}\label{Section A customized Benders decomposition method}
\begin{figure}[thbp!]
    \centering
    \begin{minipage}[t]{1\linewidth}
        \centering
        \includegraphics[width=1\textwidth]{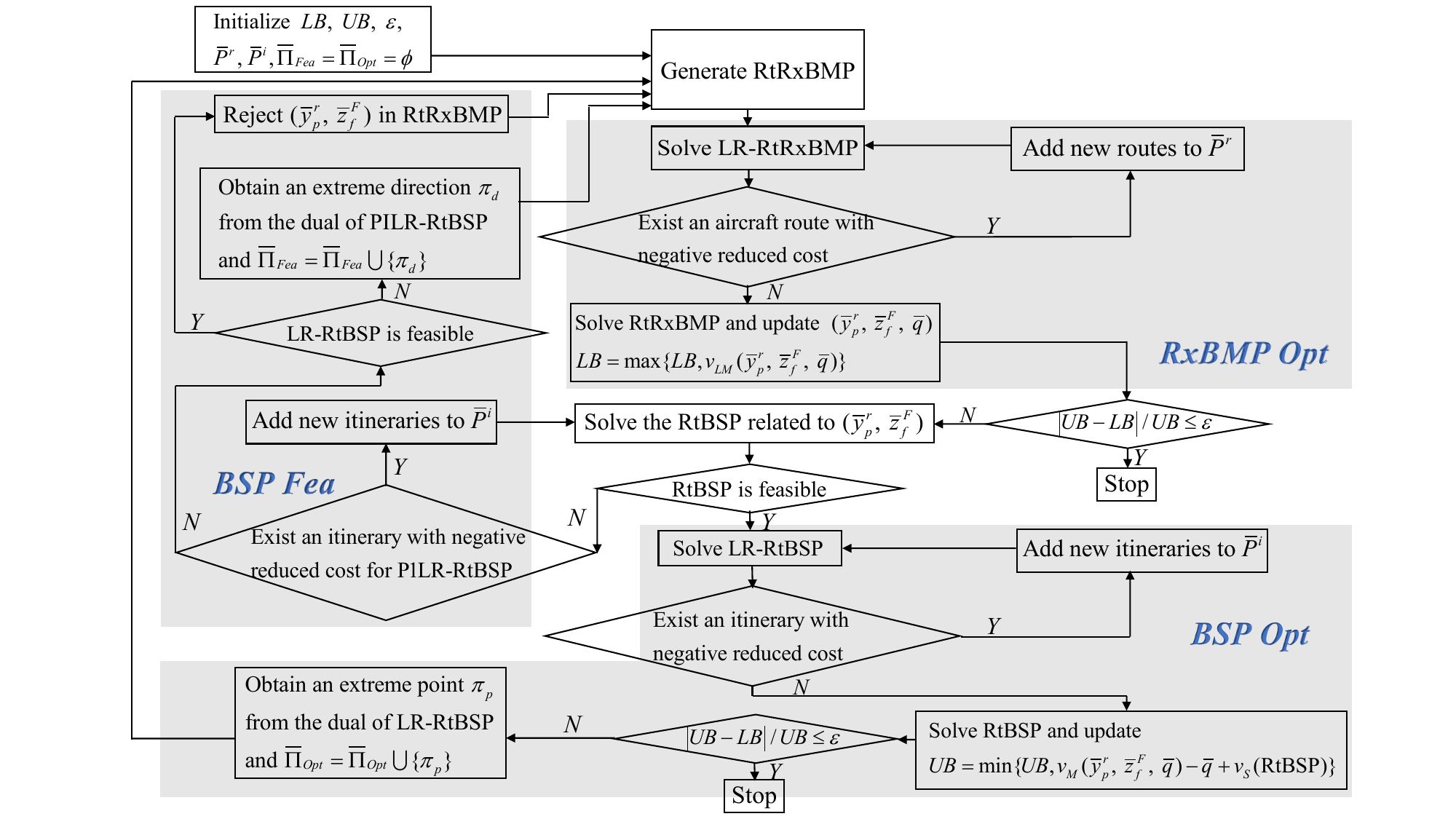}
        \caption{The flow chart of the CBD method\label{CBDflow}}
    \end{minipage}
 \end{figure}

\noindent
Having presented the Benders decomposition strategy for decomposing the IFAPRM formulation and the column generation technique for solving the Benders master problem and Benders subproblem, we now develop a customized Benders decomposition (CBD) method by combining these two effective large-scale optimization techniques. As illustrated in Figure \ref{CBDflow}, the CBD method operationalizes the proposed interactive mechanism. The CBD method mainly consists of the optimality of the RxBMP (\textit{\textbf{RxBMP Opt}}), the feasibility of the BSP (\textit{\textbf{BSP Fea}}) and the optimality of the BSP (\textit{\textbf{BSP Opt}}). Here, we focus on discussing the key points of our CBD method as follows.\par

The initialization of the sets defined in Table \ref{Notations} and the input parameters is carried out first. The initial sets for the aircraft routes in $P^r$ and the passenger itineraries in $P^i$ have a great influence on the efficiency of the CG process. Some work heuristically generates initial aircraft routes according to the planned routes or the associated shortest paths obtained by label setting algorithms (\citeauthor{maher2016}, \citeyear{maher2016}; \citeauthor{liang2018}, \citeyear{liang2018}). In this paper, we generate initial routes and itineraries by the CP solver in CPLEX, since the logic relations in the routes and itineraries can be easily formulated by CP and the resulting CP model can be effectively solved by the CP solver (\citeauthor{benoist2002}, \citeyear{benoist2002}; \citeauthor{Berbeglia2011}, \citeyear{Berbeglia2011}). For more details about CP, we refer interested readers to \citeauthor{apt2003} (\citeyear{apt2003}).

As illustrated in Figure \ref{CBDflow}, if the BSP is infeasible for a RxBMP solution $(\bar{y}_p^r, \bar{z}_f^F)$, its infeasibility can only be detected after a complete CG process. However, executing a complete CG process is time-consuming, especially for a large-scale recovery problem. To improve the solution efficiency, \citeauthor{petersen2012} (\citeyear{petersen2012}) developed an infeasibility certificate (Theorem 4.1 in their work) to detect the infeasibility of their BSP during the CG process. Once the infeasibility of their BSP is detected by this certificate, the associated CG process can be terminated immediately. Based on the structure of our BSP, we will propose a necessary and sufficient feasibility certificate in Section \ref{Section An effective feasibility certificate} that greatly reduces the runtime for checking the feasibility of the BSP, thereby accelerating the CBD method.

It is worth noting that the RxBMP and the BSP are mixed integer linear programming (MILP) problems and are solved by CG in the CBD method. When solving MILP problems by CG, some studies focus on getting an optimal solution by exact algorithms like the branch-and-price method (\citeauthor{barnhart1998}, \citeyear{barnhart1998}; \citeauthor{Levin2023}, \citeyear{Levin2023}; \citeauthor{Moreno2024}, \citeyear{Moreno2024}). Such a way outputs the best decision but usually requires substantial computational cost. To meet the real-time requirement for the integrated recovery problem, we apply the CG strategy used in \citeauthor{petersen2012} (\citeyear{petersen2012}) and \citeauthor{liang2018} (\citeyear{liang2018}) as follows. First, we generate the candidate columns by executing a complete CG process on the linear relaxation of the MILP. Then, we obtain an integer solution by solving the MILP corresponding to these columns. Such an integer solution is generally acceptable for the practical requirements, which is also applicable to our experiments.\par

\section{Acceleration Techniques}\label{Section Acceleration techniques}
\begin{figure}[thbp!]
    \centering
    \begin{minipage}[t]{1\linewidth}
        \centering
        \includegraphics[width=1\textwidth]{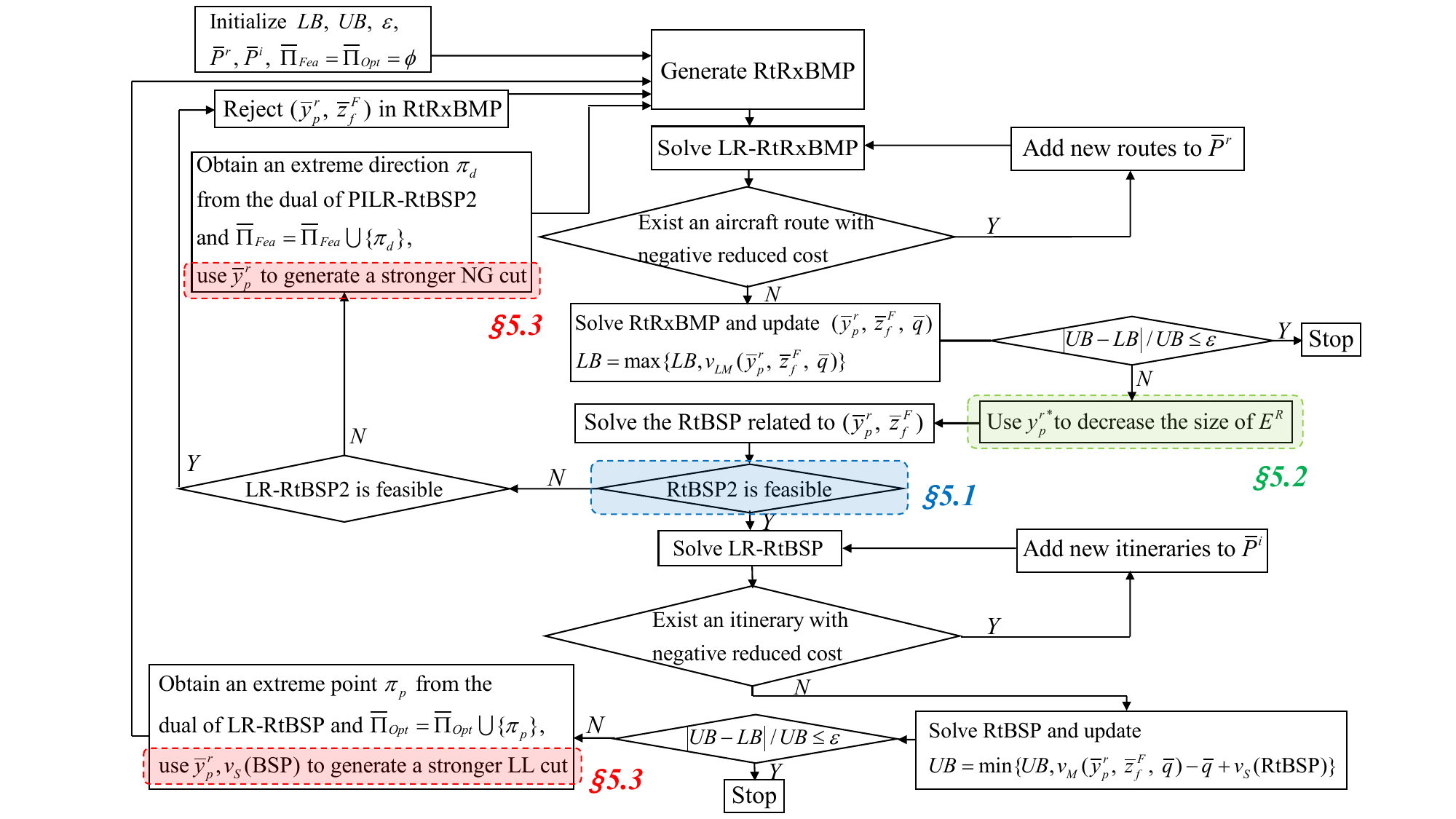}
        \caption{The flow chart of the accelerated CBD method\label{acceleratedCBDflow}}
    \end{minipage}
 \end{figure}

\noindent
Since the integrated recovery model is typically large-scale and subject to real-time requirements, we propose several acceleration techniques in this section, which are summarized in Figure \ref{acceleratedCBDflow}. The effectiveness of these acceleration techniques will be shown by the computational experiments in Section \ref{Section Computational results}.\par

\subsection{An Effective Feasibility Certificate}\label{Section An effective feasibility certificate}
\noindent
Due to the large number of elements in $P^{i}$, column generation (CG) technique is applied to solve the Benders subproblem. While CG effectively handles the large size of $P^{i}$, the infeasibility of the BSP cannot be directly detected during the CG process when the BSP is infeasible. Then, significant effort is required to execute the CG process until the infeasibility of the associated BSP can be determined. In this subsection, as summarized in Proposition \ref{BSP2feasibility}, we provide a mathematical model, denoted by BSP2, whose feasibility is equivalent to that of the BSP. Notably, the decision variables $w_p^i$ (corresponding to the elements in $P^i$) are not involved in BSP2, so we need not execute the CG process for the variables $w_p^i$ when checking the feasibility of the BSP2, which greatly improves the efficiency of checking the feasibility of our Benders subproblem.

\vspace{-0.5cm}
\begin{eqnarray}
\textbf{BSP2:}&\min &\sum_{f \in F} \sum_{f_{u} \in {FC_{f}^S}} \sum_{f_{uv} \in {FC_{f_{u}}^D}} c_{f_{uv}} x_{f_{uv}} \label{BSP2obj}\\
&\mathrm{s.t.}& \eqref{FRPRMcon1}-\eqref{FRPRMcon4} \label{BSP2con1}\\
&&
x_{f_{uv}}\in \{0,1\},\ \ \forall f \in F, f_{u} \in FC_{f}^S, f_{uv} \in FC_{f_{u}}^D. \label{BSP2con2}
\end{eqnarray}\par

\begin{proposition}
\label{BSP2feasibility}
Given a RxBMP solution $(\bar{y}_p^r, \bar{z}_f^F)$, the BSP $\eqref{FRPRMobj}-\eqref{FRPRMcon10}$ and the BSP2 $\eqref{BSP2obj}-\eqref{BSP2con2}$ have the same feasibility.
\end{proposition}

%\begin{proof}
%When the BSP is feasible, the BSP2 must be feasible since all the constraints in the BSP2 are also included in the BSP.\par
%
%Conversely, suppose $\widetilde{x}_{f_{uv}}$ is a feasible solution of the BSP2, then we can construct the following solution
%\begin{eqnarray*}(x_{f_{uv}}, z_i^{OI}, w_p^{i}) := (\widetilde{x}_{f_{uv}}, n_{i}^{OI}, \mathbf{0}),
%\end{eqnarray*}
%which is feasible for the BSP.\par
%
%Combining the above two arguments, we conclude the fact that the BSP and the BSP2 have the same feasibility.
%\end{proof}

The proof of this proposition is provided in Appendix A.1. By Proposition \ref{BSP2feasibility}, checking the feasibility of the BSP can be replaced by checking the feasibility of the BSP2. As shown in Figure \ref{acceleratedCBDflow}, when the BSP2 is feasible, we proceed to check the optimal ity of the BSP; when the BSP2 is infeasible, we generate a Benders feasibility cut associated with the BSP2 and add it to the RtRxBMP. Let $\widetilde{\Pi}_{Fea}$ denote the set of the extreme directions in the dual feasible region of the linear relaxation of the BSP2, and let
$\pi := (\pi_{f_u}^{1}, \pi_{(f_1, f_2)}^{2}, \pi_{s}^{3}, \pi_{s}^{4}) \label{the dual values of LR-BSP2}$
be the vector of the dual values related to Constraints $\eqref{FRPRMcon1}-\eqref{FRPRMcon4}$. Then the Benders feasibility cuts \eqref{Benders Fea Cut} are replaced by the following inequalities
\begin{eqnarray}
\sum_{r\in R}\sum_{p\in {P^r}}\widetilde{a}_{p}^{r}(\pi)y_p^r  + \widetilde{b}(\pi) \leq 0,
\ \ \forall \pi \in \widetilde{\Pi}_{Fea}, \label{Benders Fea Cut2}
\end{eqnarray}
where $\widetilde{a}_{p}^{r}(\pi)$ and $\widetilde{b}(\pi)$ are defined as
\begin{eqnarray*}
&&\widetilde{a}_{p}^{r}(\pi) = \sum\limits_{f \in F} \sum\limits_{f_{u} \in {FC_{f}^S}}\pi_{f_u}^{1} a_{{p,{f_u}}}^{r} + \sum\limits_{(f_1, f_2) \in {E^R}} \pi_{(f_1, f_2)}^{2}(MinTR^{r}_{(f_1,f_2)}+MaxT)b_{p,{(f_1,f_2)}}^{r},\\
&&\widetilde{b}(\pi) = -\sum\limits_{(f_1, f_2) \in {E^R}}\pi_{(f_1, f_2)}^{2}MaxT -\sum\limits_{s \in SL}(\pi_{s}^{3}u_{sdep}+\pi_{s}^{4}u_{sarr}).
\end{eqnarray*} \par

Proposition \ref{BSP2feasibility} greatly simplifies the feasibility check of the BSP. It is worth noting that the scale of the BSP2 is much smaller than that of the BSP $\eqref{FRPRMobj}-\eqref{FRPRMcon10}$. On the one hand, the variables $w_{p}^{i}$ (related to the CG process) are absent from the BSP2, which means that no CG process is required to check the feasibility of the BSP2. On the other hand, the constraint set of the BSP2 is a subset of that of the BSP. Hence, the workload on checking the feasibility is reduced, which improves the efficiency of the CBD method. The effectiveness of this feasibility certificate will be validated in Section \ref{Section Computational results}.\par

\subsection{Scale Management}\label{Section Scale management}
\noindent
When solving the BSP $\eqref{FRPRMobj}-\eqref{FRPRMcon10}$, the information in $(\bar{y}_p^r, \bar{z}_f^F)$ can be utilized to remove the redundant constraints (\citeauthor{Jiang2026}, \citeyear{Jiang2026}). From the BSP formulation, we find that the constraint in \eqref{FRPRMcon2} is redundant when the associated term $\sum\limits_{r \in R}\sum\limits_{p\in {P^r}} b_{p,{(f_1,f_2)}}^{r} \bar{y}_p^r$ is equal to zero, which means that the constraints corresponding to the non-selected flight connections in $E^{R}$ are redundant. Based on this observation, the following propositions show that most the constraints in \eqref{FRPRMcon2} are redundant when using the known information in $(\bar{y}_p^r, \bar{z}_f^F)$. We refer to this technique as \textit{scale management}.\par

\begin{proposition}
\label{Proposition Scale Management 1}
Let $n^{\bar{F}}$ be the number of flights that are not canceled. Let $\bar{R}$ denote the set of aircraft that are assigned flights, and $n^{\bar{R}}$ be its cardinality. Suppose $n^{Con}$ is the total number of flight connections in the routes of the aircraft in $\bar{R}$. Then $n^{Con}$ is equal to $n^{\bar{F}}-n^{\bar{R}}$.
\end{proposition}

%\begin{proof}
%Let $n_{r}^{\bar{R}}$ be the number of flights that are assigned to aircraft $r \in \bar{R}$. It is clear that $n^{\bar{F}} = \sum\limits_{r \in \bar{R}}n_{r}^{\bar{R}}$ holds. Then $n^{Con}$ can be calculated by summing the number of connections in each route, i.e.,
%\begin{eqnarray*}
%n^{Con} = \sum\limits_{r \in \bar{R}}(n_{r}^{\bar{R}}-1) = \sum\limits_{r \in \bar{R}}n_{r}^{\bar{R}}- \sum\limits_{r \in \bar{R}}1 = n^{\bar{F}}-n^{\bar{R}}.
%\end{eqnarray*}
%\end{proof}

The proof of Proposition \ref{Proposition Scale Management 1} is given in Appendix A.2. Based on this proposition, we provide the number of redundant constraints in \eqref{FRPRMcon2} as follows, whose proof can be found in Appendix A.3.\par

\begin{proposition}
\label{Proposition Scale Management 2}
Given a RxRMP solution $(\bar{y}_p^r, \bar{z}_f^F)$, the number of redundant constraints in Constraints \eqref{FRPRMcon2} is equal to
\begin{eqnarray*}|E^{R}|-\left(\sum_{f \in F}(1-\bar{z}_f^F) - \sum_{r\in R}\sum_{p\in {P^r}} \bar{y}_p^r\right).
\end{eqnarray*}
\end{proposition}

%\begin{proof}
%The terms $\sum\limits_{f \in F}(1-\bar{z}_f^F)$ and $\sum\limits_{r\in R}\sum\limits_{p\in {P^r}} \bar{y}_p^r$ represent the number of non-canceled flights and the number of selected aircraft, i.e., $n^{\bar{F}}$ and $n^{\bar{R}}$ in Proposition \ref{Proposition Scale Management 1}. Therefore, the number of redundant constraints in (\eqref{FRPRMcon2} is obtained by subtracting the number of selected connections from $|E^{R}|$.
%\end{proof}

We now provide a simple example to illustrate that most constraints in \eqref{FRPRMcon2} can be removed by Proposition \ref{Proposition Scale Management 2}. For each flight $f$, suppose the number of the possible successive flights for flight $f$ is $n_{f}$. Then, the number of constraints in \eqref{FRPRMcon2} is $\sum\limits_{f \in F}n_{f}$.  After applying the scale management, we only need to consider at most $|F|$ constraints in \eqref{FRPRMcon2}, which means that at least $\sum\limits_{f \in F}n_{f}-|F|$ constraints can be removed. This significantly reduces the workload of solving the BSP, thereby enhancing the efficiency of the CBD method.

\subsection{Valid Inequalities}\label{Section Valid inequalities}
\noindent
The Benders cuts \eqref{Benders Fea Cut}, \eqref{Benders Opt Cut}, and \eqref{Benders Fea Cut2} are derived from the dual information of the BSP. Since the BSP is an integer programming problem, these Benders cuts may be weak for the Benders master problem. To overcome the weakness, we will use the primal information of the BSP to generate two types of valid inequalities by exploiting the structure of IFAPRM, which are based on no-good cut and Laporte \& Louveaux cut.\par

The no-good cut (NG cut) is applicable to the problems whose linking variables are binary in Benders decomposition framework (\citeauthor{wolsey2020}, \citeyear{wolsey2020}). Let $(\bar{y}_p^r, \bar{z}_f^F)$ be a RxBMP solution that is infeasible for the BSP. The associated NG cut takes the form
\begin{eqnarray}
\sum_{(r,p)\in I_{1}(\bar{y}_p^r)}(1-{y_p^r}) + \sum_{(r,p)\in I_{0}(\bar{y}_p^r)}{y_p^r} + \sum_{f\in I_{1}(\bar{z}_f^F)}(1-z_f^F) + \sum_{f\in I_{0}(\bar{z}_f^F)}z_f^F \geq 1, \label{NG cut}
\end{eqnarray}
\noindent where $I_{1}(\bar{y}_p^r)$ and $I_{0}(\bar{y}_p^r)$ denote the sets of indexes $(r,p)$ whose associated values $\bar{y}_p^r$ are equal to 1 and 0, respectively. The definitions of $I_{1}(\bar{z}_f^F)$ and $I_{0}(\bar{z}_f^F)$ are similar to $I_{1}(\bar{y}_p^r)$ and $I_{0}(\bar{y}_p^r)$. It is noted that the NG cut only cuts off the solution $(\bar{y}_p^r, \bar{z}_f^F)$ while preserving all feasible solutions, so the NG cut is valid.\par

The Laporte \& Louveaux cut (LL cut) was originally developed for two-stage stochastic integer programming (\citeauthor{laporte1993}, \citeyear{laporte1993}) and is also applicable to problems with binary linking variables in the Benders decomposition framework. Let $(\bar{y}_p^r, \bar{z}_f^F)$ be a RxBMP solution that is feasible for the BSP. The associated LL cut is
\begin{eqnarray}
\resizebox{.9\hsize}{!}{$
q \geq q(\bar{y}_p^r, \bar{z}_f^F) + (L - q(\bar{y}_p^r, \bar{z}_f^F))\left[ \sum\limits_{(r,p)\in I_{1}(\bar{y}_p^r)}(1-{y_p^r}) + \sum\limits_{(r,p)\in I_{0}(\bar{y}_p^r)}{y_p^r} + \sum\limits_{f\in I_{1}(\bar{z}_f^F)}(1-z_f^F) + \sum\limits_{f\in I_{0}(\bar{z}_f^F)}z_f^F \right], \label{LL cut}
$}
\end{eqnarray}
\noindent where $q(\bar{y}_p^r, \bar{z}_f^F)$ is the optimal objective value of the BSP for the given solution $(\bar{y}_p^r, \bar{z}_f^F)$, and $L$ is a lower bound of $q(y_p^r, z_f^F)$ for any feasible solution $(y_p^r, z_f^F)$.  The validity of the LL cut has been discussed in \citeauthor{laporte1993} (\citeyear{laporte1993}).\par

While the NG cut \eqref{NG cut} and the LL cut \eqref{LL cut} can be directly applied to IFAPRM, the following proposition shows that these two cuts can be further strengthened based on the structure of our problem, whose proof can be found in Appendix A.4.

\begin{proposition}
\label{Proposition Stronger VI}
Given a RxRMP solution $(\bar{y}_p^r, \bar{z}_f^F)$, the following two cuts are valid and stronger than the NG cut \eqref{NG cut} and the LL cut \eqref{LL cut} respectively:
\begin{eqnarray}
&\hspace{-0.7cm}\textrm{Stronger NG cut}:&\sum_{(r,p)\in I_{1}(\bar{y}_p^r)}(1-{y_p^r}) + \sum_{(r,p)\in I_{0}(\bar{y}_p^r)}{y_p^r} \geq 1, \label{stronger NG cut}\\
&\hspace{-0.7cm}\textrm{Stronger LL cut}:&q \geq q(\bar{y}_p^r, \bar{z}_f^F) + (L - q(\bar{y}_p^r, \bar{z}_f^F))\left[ \sum_{(r,p)\in I_{1}(\bar{y}_p^r)}(1-{y_p^r}) + \sum_{(r,p)\in I_{0}(\bar{y}_p^r)}{y_p^r} \right]. \label{stronger LL cut}
\end{eqnarray}
\end{proposition}

Proposition \ref{Proposition Stronger VI} establishes that the NG cut and the LL cut can be strengthened into the valid inequalities $\eqref{stronger NG cut}$ and $\eqref{stronger LL cut}$, respectively. These tighter inequalities enhance the efficiency of solving the Benders master problem. Their effectiveness is verified through computational experiments in Section \ref{Section Computational results}.

\section{Numerical Experiments}\label{Section Numerical experiments}
\noindent
We perform a series of experiments to demonstrate the effectiveness of the proposed methodology. Our experiments are based on the real-world data from the ROADEF 2009 challenge (\url{https://www.roadef.org/challenge/2009/en/instances.php}). These experiments are run on a laptop with a 2.60 GHz Intel i7 CPU and implemented in Java. In addition, we call CPLEX 12.9 to solve the involved linear programmings, mixed integer linear programmings, and constraint programmings.\par

\subsection{Data Description}\label{Section Data description}
\noindent
Ten instances with different scales and disruption scenarios used in our experiments are summarized in Table \ref{Characteristics of data set}. The recovery time window is set to be one day.  Maintenance tasks are regarded as the special flights whose departure times and arrival times cannot be changed. To reduce the scale of passenger recovery, passengers are grouped into several sets according to their original itineraries (flight sequences). If two passengers share the same original itinerary, they are placed in the same set, though their actual itineraries may differ in a recovery plan.\par

%\vspace{-1em}
\begin{table}[!htbp]
  \centering
  \renewcommand{\baselinestretch}{1.2}
  \caption{Characteristics of data set in computational experiments\label{Characteristics of data set}}
  \small
    {\resizebox{1\linewidth}{!}{
    \begin{tabular}{cccccccc}
    \toprule
    Instance & No. of flights & No. of aircraft & No. of FTs\footnotemark[1] & No. of airports & No. of PAXs\footnotemark[2] & No. of OIs\footnotemark[3] & Disruption \\
    \midrule
    1     & 12    & 3     & 1     & 7     & 656   & 12    & FCon\footnotemark[4] \\
    2     & 39    & 7     & 2     & 14    & 884   & 44    & FCon \\
    3     & 73    & 13    & 3     & 19    & 2,508  & 99    & FDD\footnotemark[5] \\
    4     & 121   & 21    & 4     & 22    & 5,332  & 174   & FCon \\
    5     & 121   & 21    & 4     & 22    & 5,332  & 174   & FDD \\
    6     & 121   & 21    & 4     & 22    & 5,332  & 174   & FCan\footnotemark[6] \\
    7     & 177   & 31    & 4     & 25    & 9,616  & 263   & FCan \\
    8     & 272   & 45    & 5     & 28    & 19,704 & 532   & FCan \\
    9     & 372   & 63    & 10    & 34    & 27,189 & 752   & FDD \\
    10    & 372   & 63    & 10    & 34    & 27,189 & 752   & FCan \\
    \bottomrule
    \multicolumn{8}{l}{\footnotemark[1]FTs: fleet types. \footnotemark[2]OIs: original itineraries.  \footnotemark[3]PAXs: passengers.}\\
    \multicolumn{8}{l}{\footnotemark[4]FCon: flow control. \footnotemark[5]FDD: flight departure delay. \footnotemark[6]FCan: flight cancellation.}\\
    \end{tabular}
    }}
\end{table}
%\vspace{-1em}

%\vspace{0pt}
\begin{table}[!htbp]
  \centering
  \renewcommand{\baselinestretch}{1.1}
  \caption{Benchmark parameters used in experiments\label{Benchmark parameters used in experiments}}
  \small
   \setlength{\tabcolsep}{5mm}{
    \begin{tabular}{l|c}
    \toprule
    Maximum allowable delay for each flight (min) & 120 \\
    Length of sparse copy interval (min) & 30 \\
    Length of dense copy interval (min) & 5 \\
    Cost of cancellation per flight (\$) & 25,000 \\
    Cost of flight delay (\$/min) & 100 \\
    Cost of swap per tail assignment (\$) & 0 \\
    Non-cruise time for each flight (min) & 30\\
    Cost of fuel (\$/kg) & 1 \\
    CO$_2$ emission constant $\kappa$ in the functions \eqref{Fuel Formulation} & 3.15 \\
    Cost of carbon emission (\$/kg) & 0.02 \\
    Cost of itinerary cancellation per passenger (\$) & 2,500 \\
    Cost of passenger delay (\$/min) & 0.64 \\
    Cost of changing passenger itineraries (\$) & 0 \\
    Maximum number of flights in each passenger itinerary & 4 \\
    \bottomrule
    \end{tabular}
    }
\end{table}

The key parameters of our experiments are shown in Table \ref{Benchmark parameters used in experiments}. The cost parameters of flight cancellation, flight swap, passenger itinerary cancellation, passenger delay, and itinerary change are the same as those in \citeauthor{petersen2012} (\citeyear{petersen2012}). The cost of flight delay is obtained from \citeauthor{ball2010} (\citeyear{ball2010}). These cost parameters depend on the preferences of airline companies and other factors. To show their effects on the recovery performance, we conduct the sensitivity analyses in Section \ref{Section Sensitivity analysis results}.\par

The cost of controlling cruise speeds is calculated by the functions in \citeauthor{akturk2014} (\citeyear{akturk2014}), i.e., the fuel consumption function $F(v)$, the cost function for the change of fuel consumption $\Delta \textrm{Fuel Cost} (v)$, and the cost function for the change in carbon emission $\Delta \textrm{Carbon Emission Cost} (v)$:

%\vspace*{0em}
\begin{table*}[h]
\begin{numcases}{}\label{Fuel Formulation}
F(v) = d^{cr}\cdot(c_{1}v^{2}+c_{2}v+\frac{c_{3}}{v^{2}}+\frac{c_{4}}{v^{3}}), \notag\\
\Delta \textrm{Fuel Cost} (v) = c_{\textrm{fuel}}\cdot\left(F(v)-F(v^{0})\right), \\
\Delta \textrm{Carbon Emission Cost} (v) = c_{\textrm{CO$_2$}}\cdot\kappa\cdot\left(F(v)-F(v^{0})\right), \notag
\end{numcases}
\end{table*}
\noindent
where $v$ denotes the cruise speed of a flight. Coefficients $c_1, ..., c_4 > 0$ are influenced by the specific characteristics of the aircraft and the cruise environment. In addition, $d^{cr}$ and $v^{0}$ denote the distance at the cruise stage and the planned cruise speed, respectively. The above coefficients are also obtained from  \citeauthor{akturk2014} (\citeyear{akturk2014}). In this paper, the cruise speed control is considered in a discrete way. For each flight, the candidate cruise speeds are uniformly chosen from the interval of cruise speeds. Specifically, suppose the interval is $[v^{0}, \alpha\cdot v^{0}]$, then the set of candidate cruise speeds is
\begin{eqnarray}
\left\{v~|~v=v^{0}+k\cdot\frac{\alpha \cdot v^{0}-v^{0}}{\max\{K-1,1\}},~k=0,..,K-1\right\}.  \label{the set of candidate cruise speeds}
\end{eqnarray}

\noindent
In our experiments, the coefficient related to the maximum allowable cruise speed is set to 1.1, as in \citeauthor{akturk2014} (\citeyear{akturk2014}), and the number of candidate cruise speeds $K$ is set to 5. We will discuss the effects of the parameters $\alpha$ and $K$ in Section \ref{Section The effect of cruise speed control}. Accordingly, given a cruise speed from the set \eqref{the set of candidate cruise speeds}, the cost related to the cruise speed control is calculated by the equations in \eqref{Fuel Formulation}.\par

\subsection{Computational Results}\label{Section Computational results}
\noindent
To compare the performance of our interactive mechanism with that of the conventional flight copy approach, we formulate an IFAPRM with all dense flight copies (IFAPRM-AD) as follows. The objective function and constraints of IFAPRM-AD are similar to those of IFAPRM except that IFAPRM-AD uses only dense flight copies.\par

%\vspace*{-1em}
\begin{eqnarray}
\nonumber
&&\hspace{-1.5cm}\textbf{IFAPRM-AD formulation:}\\
\nonumber
&\min &\sum_{f\in F}(d_{f}^{F} - c_{f}^{OF})z_{f}^{F} + \sum_{r\in R}\sum_{p\in {P^r}}c_p^ry_p^r +
\sum_{f \in F} \sum_{f_{uv} \in FC_{f}^{D}} c_{f_{uv}} x_{f_{uv}} \label{M-IFAPRMobj}\\
&&
 + \sum_{i \in OI}d_{i}^{OI}z_{i}^{OI}+\sum_{i \in OI}\sum_{p \in P^i}c_{i}^{DT}DT_{p}^{i}w_{p}^{i} + \sum_{i \in OI}\sum_{p \in P^i}c_{i}^{C}n_{i,p}^{C}w_{p}^{i}\\
&\mathrm{s.t.} &\sum_{f_{uv} \in FC_{f}^{D}}x_{f_{uv}}+z_{f}^{F}=1,\ \ \forall f\in {F},\label{M-IFAPRMcon1}\\
&&
\sum_{p\in {P^r}} y_p^r \leq 1,\ \ \forall r\in R,\label{M-IFAPRMcon2}\\
&&
x_{f_{uv}} = \sum_{r\in R}\sum_{p\in {P^r}} a_{{p,{f_{uv}}}}^{r} y_p^r,\ \ \forall f\in {F}, f_{uv} \in FC_{f}^{D},\label{M-IFAPRMcon3}\\
&&
\sum_{f \in F} \sum_{f_{uv} \in FC_{f}^{D}} a_{f_{uv}}^{sdep} x_{f_{uv}} \leq u_{sdep},\ \ \forall s\in {SL},\label{M-IFAPRMcon4}\\
&&
\sum_{f \in F} \sum_{f_{uv} \in FC_{f}^{D}} a_{f_{uv}}^{sarr} x_{f_{uv}} \leq u_{sarr},\ \ \forall s\in {SL},\label{M-IFAPRMcon5}\\
&&
\sum_{p\in {P^i}} w_p^i + z_i^{OI} = n_{i}^{OI},\ \ \forall i\in OI,\label{M-IFAPRMcon6}\\
&&
\sum_{i \in OI}\sum_{p \in P^i} a_{{p,{f_{uv}}}}^{i} w_p^i \leq MaxCap \cdot x_{f_{uv}},\ \ \forall f \in F, f_{uv} \in {FC_{f}^D},\label{M-IFAPRMcon7}\\
\nonumber
&&
\sum_{i \in OI}\sum_{p \in P^i}\sum_{f_{uv} \in {FC_{f}^D}} a_{{p,{f_{uv}}}}^{i}w_{p}^{i}
\leq \\
&&
\hspace{3cm} \sum_{r \in R}Cap_r \sum_{p \in P^r}\sum_{f_{uv} \in FC_{f}^D}a_{{p,{f_{uv}}}}^{r} \bar{y}_p^r + MaxCap\cdot \bar{z}_{f}^{F},\ \ \forall f \in F,\label{M-IFAPRMcon8}\\
&&
z_{f}^{F}\in \{0,1\},\ \ \forall j\in {F},\label{M-IFAPRMcon9}\\
&&
y_p^r\in \{0,1\},\ \ \forall r\in R, p\in P^r,\label{M-IFAPRMcon10}\\
&&
x_{f_{uv}}\in \{0,1\},\ \ \forall f \in F, f_{uv} \in FC_{f}^{D},\label{M-IFAPRMcon11}\\
&&
z_i^{OI} \in \mathbb{Z_{+}},\ \ \forall i\in OI,\label{M-IFAPRMcon12}\\
&&
w_p^{i} \in \mathbb{Z_{+}},\ \ \forall i\in OI, p \in P^i.\label{M-IFAPRMcon13}
\end{eqnarray}\par

We carry out experiments using the CBD method and CPLEX to solve IFAPRM and IFAPRM-AD, yielding four comparative  variants: SD-CBD, AD-CBD, SD-CPLEX, and AD-CPLEX. Their results are reported in Table \ref{Results of CBD and CPLEX}. When solving IFAPRM-AD by the CBD method (AD-CBD), the first three terms of objective function \eqref{M-IFAPRMobj}, Constraints $\eqref{M-IFAPRMcon1}-\eqref{M-IFAPRMcon5}$ and $\eqref{M-IFAPRMcon9}-\eqref{M-IFAPRMcon11}$ are put in the Benders master problem, and the remaining terms are put in the Benders subproblem. When using CPLEX to solve IFAPRM and IFAPRM-AD (SD-CPLEX and AD-CPLEX, respectively), aircraft routes and passenger itineraries are generated as follows. Besides the initial routes and itineraries, CPLEX finds as many additional routes and itineraries as possible within 5 seconds for each aircraft and 2 seconds for each passenger. In addition, to meet practical requirements, the maximum allowable runtime of four variants is 30 minutes, which is consistent with \citeauthor{petersen2012} (\citeyear{petersen2012}). The tolerance of the optimality gap $(UB-LB)/UB\times100\%$ for the CBD method is set as 5\%.\par

\begin{table}[!htbp]
  \centering
  \renewcommand{\arraystretch}{1.15}
  \caption{Results of the CBD method and CPLEX for solving the IFAPRM and IFAPRM-AD \label{Results of CBD and CPLEX}}
    {\resizebox{\linewidth}{!}{
   \begin{tabular}{l|cccccccccc}
    \toprule
    Instance & 1     & 2     & 3     & 4     & 5     & 6     & 7     & 8     & 9     & 10 \\
    \midrule
    \multicolumn{11}{l}{Solution metrics by SD-CBD\footnotemark[1]}\\
    \midrule
    No. of canceled flights & 2     & 2     & 3     & 0     & 3     & 3     & 6     & 4     & 0     & 4 \\
    Flight delays (min) & 235   & 160   & 49    & 235   & 24    & 5     & 30    & 0     & 321   & 0 \\
    No. of unassigned passengers & 40    & 0     & 35    & 0     & 35    & 32    & 0     & 0     & 0     & 37 \\
    Passenger delays (min) & 43,150 & 16,175 & 5,274  & 10,810 & 5,124  & 4,965  & 61,330 & 20,370 & 23,366 & 35,800 \\
    Change of used fuel (kg) & -3,028  & -1,606  & -2,533  & 166   & -2,558  & -2,586  & -12,898  & -8,205  & 427   & -10,496  \\
    Change of carbon emission (kg) & -9,538  & -5,059  & -7,979  & 522   & -8,059  & -8,147  & -40,629  & -25,846  & 1,344  & -33,062  \\
    Recovery cost (\$) & 197,897  & 74,645  & 168,083  & 30,595  & 165,460  & 155,928  & 178,541  & 104,314  & 47,508  & 204,255  \\
    Optimality gap (\%) & 4.80  & 4.63  & 2.86  & 1.47  & 1.56  & 3.44  & 3.37  & 2.10  & 3.25  & 4.66  \\
    CPU time (s) & 2     & 11    & 14    & 59   & 56   & 47   & 295   & 490  & 572  & 548  \\
    \midrule
    \multicolumn{11}{l}{Solution metrics by AD-CBD\footnotemark[2]}\\
    \midrule
    No. of canceled flights & 2     & 2     & 3     & 0     & 3     & 3     & 6     & 4     & 0     & 4 \\
    Flight delays (min) & 229   & 147   & 54    & 235   & 29    & 5     & 0     & 0     & 319   & 0 \\
    No. of unassigned passengers & 40    & 0     & 33    & 0     & 35    & 33    & 0     & 0     & 0     & 33 \\
    Passenger delays (min) & 42,814 & 15,876 & 5,670  & 10,810 & 5,254  & 5,240  & 60,100 & 17,200 & 22,498 & 36,659 \\
    Change of used fuel (kg) & -2,966  & -1,497  & -2,507  & 170   & -2,501  & -2,593  & -12,898 & -8,205 & 388   & -10,488  \\
    Change of carbon emission (kg) & -9,343  & -4,714  & -7,896  & 536   & -7,877  & -8,169  & -40,629 & -25,846  & 1,221  & -33,038  \\
    Recovery cost (\$) & 197,148  & 73,270  & 163,864  & 30,599  & 166,104  & 158,597  & 174,753 & 102,286  & 46,711  & 194,813  \\
    Optimality gap (\%) & 0.13  & 3.02  & 0.62  & 0.18  & 0.59  & 4.14  & 4.31  & 0.50  & 2.65  & 0.19  \\
    CPU time (s) & 13    & 155   & 294   & 465  & 349  & 606  & 823  & 1,176  & 1,800+  & 1,800+  \\
    \midrule
    \multicolumn{11}{l}{Solution metrics by SD-CPLEX\footnotemark[3]}\\
    \midrule
    No. of canceled flights & 2     & 2     & 3     & 0     & 3     & 3     & n/a\footnotemark[6]   & n/a   & n/a   & n/a \\
    Flight delays (min) & 232   & 147   & 48    & 235   & 24    & 5     & n/a   & n/a   & n/a   & n/a \\
    No. of unassigned passengers & 40    & 0     & 35    & 0     & 35    & 32    & n/a   & n/a   & n/a   & n/a \\
    Passenger delays (min) & 42,937  & 15,876  & 5,268  & 10,810 & 5,124  & 4,965  & n/a   & n/a   & n/a   & n/a \\
    Change of used fuel (kg) & -2,997  & -1,495  & -2,529  & 124   & -2,579  & -2,586  & n/a   & n/a   & n/a   & n/a \\
    Change of carbon emission (kg) & -9,441  & -4,710  & -7,967  & 389   & -8,124  & -8,147  & n/a   & n/a   & n/a   & n/a \\
    Recovery cost (\$) & 197,494  & 73,271  & 167,983  & 30,550 & 165,438  & 155,928  & n/a   & n/a   & n/a   & n/a \\
    CPU time (s) & 94   & 496  & 843  & 1,800+\footnotemark[5]  & 1,428  & 1,800+  & --\footnotemark[7] & -- & --     & -- \\
    \midrule
    \multicolumn{11}{l}{Solution metrics by AD-CPLEX\footnotemark[4]}\\
    \midrule
    No. of canceled flights & 2     & 2     & 3     & 0     & 3     & 3     & 6     & 4     & n/a     & n/a \\
    Flight delays (min) & 229   & 147   & 54    & 235   & 24    & 5     & 0     & 0     & n/a     & n/a \\
    No. of unassigned passengers & 40    & 0     & 33    & 0     & 35    & 33    & 0     & 0     & n/a     & n/a \\
    Passenger delays (min) & 42,569  & 15,876  & 5,670  & 10,810 & 5,124  & 5,240  & 59,916 & 16,970 & n/a     & n/a \\
    Change of used fuel (kg) & -2,925  & -1,506  & -2,507  & 170   & -2,579  & -2,593  & -12,824 & -8,165  & n/a     & n/a  \\
    Change of carbon emission (kg) & -9,212  & -4,745  & -7,896  & 536   & -8,124  & -8,169  & -40,395 & -25,720  & n/a     & n/a  \\
    Recovery cost (\$) & 197,035  & 73,260  & 163,864  & 30,599  & 165,438  & 158,597  & 174,714 & 102,181  & n/a     & n/a  \\
    CPU time (s) & 71   & 129   & 172   & 190   & 193   & 218   & 446  & 587  & --  & --  \\
    \bottomrule
    \multicolumn{11}{l}{\footnotemark[1]SD-CBD: IFAPRM with the CBD method.~~\footnotemark[2]AD-CBD: IFAPRM-AD with the CBD method.}\\
    \multicolumn{11}{l}{\footnotemark[3]SD-CPLEX: IFAPRM with CPLEX.~~\footnotemark[4]AD-CPLEX: IFAPRM-AD with CPLEX.}\\
    \multicolumn{11}{l}{\footnotemark[5]1,800+: out of the maximum allowable runtime, 1,800 seconds.~~\footnotemark[6]n/a: not applicable.~~\footnotemark[7]--: out of memory.}
    \end{tabular}
    }}
\end{table}

\subsubsection{Insights from Computational Results.}
\noindent
We obtain the following observations from Table \ref{Results of CBD and CPLEX}. First, all test instances are solved by SD-CBD within 10 minutes under the given optimality gap restriction. Even for the largest-scale instances (\#9 and \#10), SD-CBD completes within 572 and 548 seconds, respectively. This demonstrates that the proposed methodology---the sparse-dense flight copy-based mechanism, along with its corresponding modeling method, algorithm, and acceleration techniques---meets the practical requirement of addressing the integrated recovery problem under time pressure.

Second, by comparing the results of SD-CBD with those of AD-CBD, we observe that SD-CBD achieves broadly comparable solution quality in terms of recovery cost and optimality gap, while requiring significantly less runtime. For example, for Instance \#4, the recovery costs are \$30,595 (SD-CBD) and \$30,599 (AD-CBD), with gaps of 1.47\% and 0.18\%, respectively, yet SD-CBD runs in only 59 seconds compared to 465 seconds for AD-CBD. This advantage becomes even more pronounced for large-scale instances. For Instances \#9 and \#10, AD-CBD exceeds the maximum allowable runtime of 1,800 seconds, rendering it impractical for real-time recovery. In contrast, SD-CBD solves both instances within 10 minutes, demonstrating its suitability for practical disruption operations. The inferior computational efficiency of AD-CBD is attributable to the model itself rather than the solution method: since aircraft routes in IFAPRM-AD consist entirely of dense flight copies, the number of possible routes grows explosively, dramatically increasing the computational burden. In contrast, the sparse-dense flight copy approach in IFAPRM keeps both the sparse and dense flight copy networks compact, enabling the Benders decomposition framework to operate efficiently.

Third, from a comparison of SD-CBD and SD-CPLEX, we find that the CBD method is more scalable and robust for solving the integrated recovery problem than CPLEX. For small-scale to medium-scale instances (\#1$-$\#6), SD-CBD and SD-CPLEX yield comparable recovery costs, but SD-CBD is consistently faster (e.g., Instance \#3: 14s vs. 843s; Instance \#4: 59s vs. 1,800+s). More importantly, for larger-scale instances (\#7$-$\#10), SD-CPLEX fails to find any feasible solution due to memory limitations (denoted as ``n/a"), whereas SD-CBD successfully solves all of them. This robustness is attributed to the decomposition strategy: the CBD method divides IFAPRM into several subproblems whose scales are much smaller than that of the full model, thus requiring less memory and computational effort.

Fourth, comparing the results of SD-CBD with those of AD-CPLEX (the conventional approach) reveals the combined benefit of the sparse-dense flight copy approach and the customized solution method. For small-scale instances (\#1$-$\#3), AD-CPLEX produces slightly better recovery costs than SD-CBD, but SD-CBD runs considerably faster. For medium-scale instances with 121 flights  (\#4$-$\#6), SD-CBD  reduces runtime by over 68\% and yields solutions with comparable quality. For the largest-scale instances with 372 flights (\#9 and \#10), AD-CPLEX runs out of memory and cannot produce any feasible solution, while SD-CBD delivers feasible recovery plans within 10 minutes. This demonstrates that the proposed interactive mechanism strikes an effective balance between solution quality and computational tractability, making it well-suited for large-scale integrated recovery problems under real-time requirements.

The above observations collectively demonstrate that the proposed interactive mechanism, operationalized through the IFAPRM formulation, the CBD method, and acceleration techniques, delivers high-quality recovery solutions within practical runtime across all test instances, confirming its suitability for real-world airline disruption management.

\subsubsection{The Effect of Cruise Speed Control.}\label{Section The effect of cruise speed control}
\begin{figure}[!htbp]
   \hspace{-0.15cm}
   \begin{minipage}{0.49\linewidth}
    \centering
    \includegraphics[width=1\textwidth]{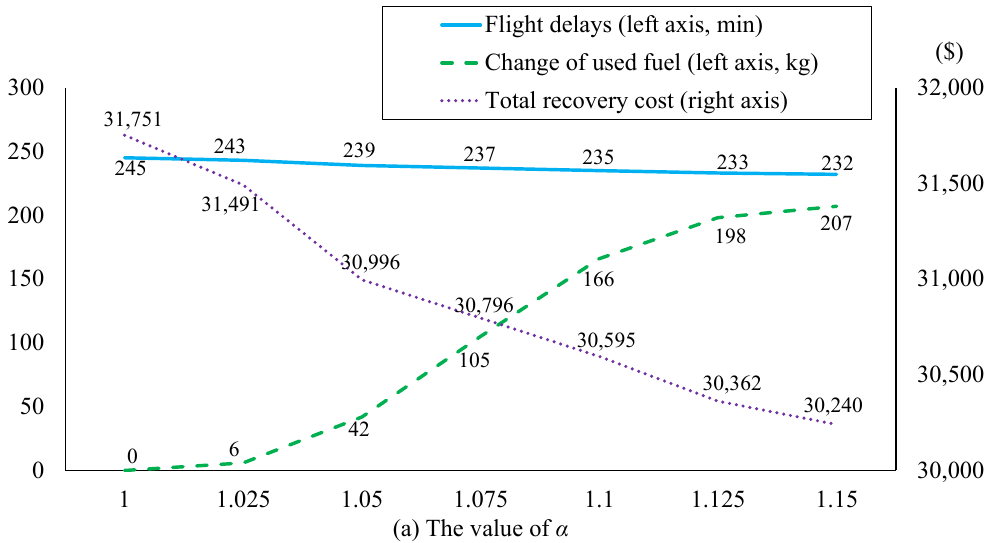}
   \end{minipage}
   \qquad
   \hspace{-0.8cm}
   \begin{minipage}{0.5\linewidth}
    \centering
    \vspace{-0.10cm}
    \includegraphics[width=1\textwidth]{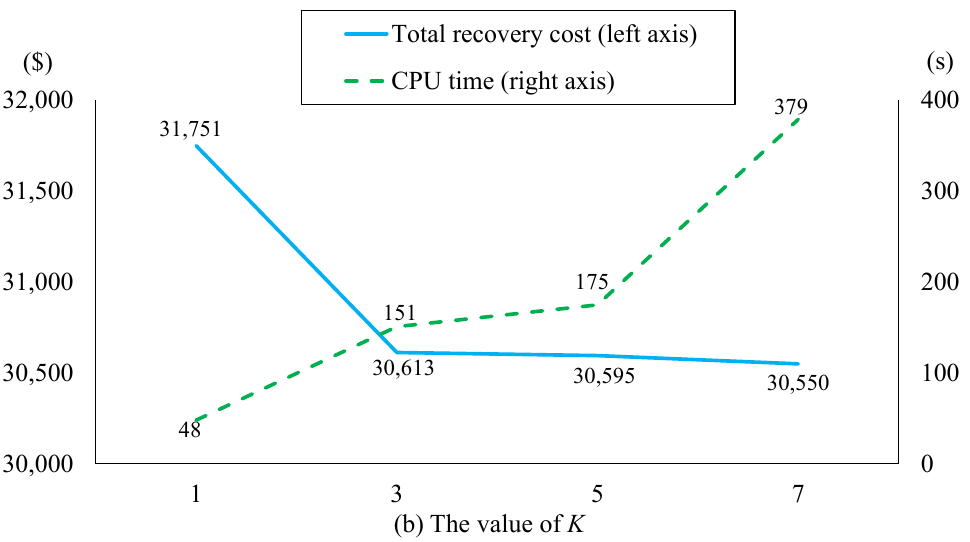}
   \end{minipage}
   \caption{Sensitivity analysis for the values of $\alpha$ and $K$ in the set \eqref{the set of candidate cruise speeds} \label{value_alphavalue_k}}
\end{figure}

\noindent
Cruise speed control is an effective recovery option that can enhance both the efficiency and flexibility of the recovery plan by compressing flying times and providing more adjustable departure/arrival opportunities. In our integrated recovery model IFAPRM, cruise speed control is incorporated as one of the recovery options. To demonstrate its effect on recovery performance, we vary the values of $\alpha$ and $K$ in the set \eqref{the set of candidate cruise speeds}. The associated results for Instance \#4, a medium-scale instance in our experiments, are illustrated in Figure \ref{value_alphavalue_k}. Other instances exhibit similar behaviors.

As shown in Figure \ref{value_alphavalue_k}(a), comparing the case $\alpha = 1$ (i.e., no cruise speed control) with other values, we find that cruise speed control can indeed improve recovery performance. As the maximum allowable cruise speed increases (by increasing $\alpha$), the total delay time gradually decreases, as expected, despite some increase in fuel consumption, leading to a lower total recovery cost.

The following observations are obtained from Figure \ref{value_alphavalue_k}(b). Similar to the value $\alpha = 1$ in Figure \ref{value_alphavalue_k}(a), by comparing the value $K = 1$ with other values, we further find that controlling cruise speeds can reduce the total recovery cost. Meanwhile, as $K$ increases, the required runtime grows significantly. The underlying reasons are as follows. Introducing more candidate cruise speeds provides more feasible recovery plans, thus leading to a lower total recovery cost. However, the scale of IFAPRM also increases, requiring more computational effort to solve the problem. Thus, an appropriate number of candidate cruise speeds needs to be chosen to balance the trade-off between real-time computational feasibility and recovery plan quality.

\subsubsection{The Effectiveness of Acceleration Techniques.}
\noindent
To demonstrate the effectiveness of the acceleration techniques proposed in Section \ref{Section Acceleration techniques}, we conduct experiments using different combinations of acceleration techniques. Table \ref{acceleration techniques experiments} provides the associated results for Instance \#4 in Table \ref{Characteristics of data set}, and similar results are observed for other instances.

\begin{table}[!htbp]
  \centering
  \renewcommand{\baselinestretch}{1.2}
  \caption{The results of different acceleration techniques for the CBD method \label{acceleration techniques experiments}}
  \small
  \setlength{\tabcolsep}{2mm}{
    {\resizebox{1\linewidth}{!}{
    \begin{tabular}{p{0.33\textwidth}|p{2.8cm}<{\centering}p{3.2cm}<{\centering}p{2.5cm}<{\centering}p{1.5cm}<{\centering}p{1.5cm}<{\centering}p{3.5cm}<{\centering}p{1.5cm}<{\centering}}
    \toprule
              & A     & B     & C$_1$     & C$_2$     & C$_3$     & C$_4$    & D \\
    \midrule
    Situation & No Fea. Cert.\footnotemark[1] & No Scale Manag.\footnotemark[2] & BC\footnotemark[3] & VI\footnotemark[4] & SVI\footnotemark[5] & BC + VI\footnotemark[6] & All\footnotemark[7]  \\
    \midrule
   No. of canceled flights & 0     & 0     & 0     & n/a\footnotemark[8]   & n/a   & 0     & 0 \\
    Flight delays (min) & 235   & 235   & 235   & n/a   & n/a   & 235   & 235 \\
    No. of unassigned passengers & 0     & 0     & 0     & n/a   & n/a   & 0     & 0 \\
    Passenger delays (min) & 10,810 & 10,810 & 10,810 & n/a   & n/a   & 10,810 & 10,810 \\
    Change of used fuel (kg) & 166   & 166   & 177   & n/a   & n/a   & 166   & 166  \\
    Change of carbon emission (kg) & 522   & 522   & 556   & n/a   & n/a   & 522   & 522  \\
    Recovery cost (\$) & 30,595  & 30,595  & 30,606  & n/a & n/a & 30,595  & 30,595  \\
    Optimality gap (\%) & 1.47  & 1.47  & 2.05  & n/a   & n/a   & 1.36  & 1.47  \\
    CPU time (s) & 112   & 81   & 83   & 1,800+\footnotemark[9]  & 1,800+  & 76   & 59  \\
    \bottomrule
    \multicolumn{8}{l}{\footnotemark[1]No Fea. Cert.: the situation using all acceleration techniques except the feasibility certificate proposed in Section \ref{Section An effective feasibility certificate}.}\\
    \multicolumn{8}{l}{\footnotemark[2]No Scale Manag.: the situation using all acceleration techniques except the scale management proposed in Section \ref{Section Scale management}.}\\
    \multicolumn{8}{l}{\footnotemark[3]BC: the situation using the feasibility certificate, the scale management, and the Benders cuts \eqref{Benders Opt Cut} and \eqref{Benders Fea Cut2}.}\\
    \multicolumn{8}{l}{\footnotemark[4]VI: the situation using the feasibility certificate, the scale management, and the valid inequalities \eqref{NG cut} and \eqref{LL cut}. }\\
    \multicolumn{8}{l}{\footnotemark[5]SVI: the situation using the feasibility certificate, the scale management, and the stronger valid inequalities \eqref{stronger NG cut} and \eqref{stronger LL cut}.}\\
    \multicolumn{8}{l}{\footnotemark[6]BC + VI: the situation using the situation using the feasibility certificate, the scale management, the Benders cuts, and the valid inequalities.}\\
    \multicolumn{8}{l}{\footnotemark[7]All: the situation using the situation using the feasibility certificate, the scale management, the Benders cuts, and the stronger valid inequalities.}\\
    \multicolumn{8}{l}{\footnotemark[8]n/a: fail to find a feasible solution.~~\footnotemark[9]1,800+: out of the maximum allowable runtime of solving IFAPRM by the CBD method, 1,800 seconds.}
    \end{tabular}
    }}
    }
\end{table}

The comparison of Column A and Column D shows that the proposed feasibility certificate can save 47\% of the runtime. This substantial reduction is because we only need to check the feasibility of BSP2 $\eqref{BSP2obj}-\eqref{BSP2con2}$ instead of the full BSP. The scale of BSP2 is much smaller than that of the BSP. On the one hand, BSP2 contains only the constraints that affect the feasibility of the BSP. On the other hand, BSP2 has far fewer variables since it does not include the variables $w_{p}^{i}$ (related to the CG process), meaning that the CG process is not required when checking the feasibility of BSP2. Due to these two factors, the runtime for checking the feasibility of the Benders subproblem in Column D is much shorter, leading to a significant reduction in the total runtime.

Compared with Column B, the runtime of Column D decreases by approximately 27\% with the help of scale management. By applying this acceleration technique, many redundant constraints in Constraints \eqref{FRPRMcon2} are removed for Column D, thereby reducing the workload of solving the BSP.

From Columns C$_1-$C$_4$ and Column D, the following observations are obtained regarding the cuts used in the accelerated CBD method. First, from the result of Column C$_1$, we find that using only Benders cuts (without using the proposed valid inequalities) is not so effective, which is attributed to the fact that Column C$_1$ only uses the dual information of the BSP. Compared with Column C$_1$, the runtime of Column D is reduced by roughly 29\%. Second, Columns C$_2$ and C$_3$ show that a feasible solution cannot be found within the maximum allowable runtime when using the valid inequalities $\eqref{NG cut}-\eqref{LL cut}$ or the stronger valid inequalities $\eqref{stronger NG cut}-\eqref{stronger LL cut}$ instead of Benders cuts. This is because the cuts $\eqref{NG cut}$ and $\eqref{stronger NG cut}$ only use the primal information of the BSP to eliminate the associated infeasible solution $(\bar{y}_p^r, \bar{z}_f^F)$ but do not restrict other infeasible solutions. These two observations imply that using only dual or primal information leads to poor performance. As shown in Columns C$_4$ and D, this performance can be improved by simultaneously utilizing both dual and primal information. Third, from the comparison of Columns C$_4$ and D, we find that the runtime of Column D can be reduced by using the stronger valid inequalities $\eqref{stronger NG cut}-\eqref{stronger LL cut}$. This improvement is attributable to the stronger valid inequalities derived from sufficiently exploiting the structure of IFAPRM.

\subsection{Sensitivity Analysis Results}\label{Section Sensitivity analysis results}
The cost parameters in IFAPRM are related to airline preferences and resource availability constraints. In this subsection, we conduct sensitivity analyses on the cost parameters to explore their effects on recovery performance. The associated experiments are also carried out on Instance \#4, and the other instances exhibit similar performance.

\subsubsection{Sensitivity Analysis for the Unit Delay Cost and the Unit Cancellation Cost.}
\begin{figure}[!htbp]
   \hspace{-0.1cm}
   \begin{minipage}{0.5\linewidth}
    \centering
    \includegraphics[width=1\textwidth]{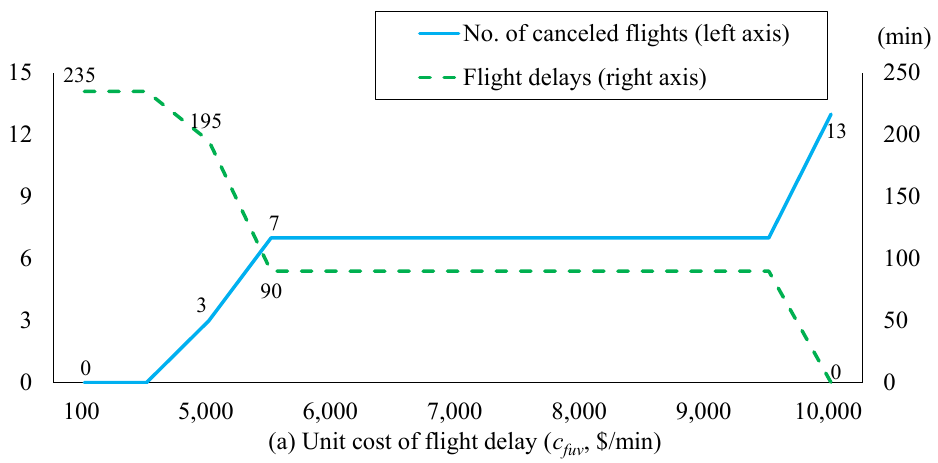}
   \end{minipage}
   \qquad
   \hspace{-0.8cm}
   \begin{minipage}{0.48\linewidth}
    \centering
    \vspace{-0.135cm}
    \includegraphics[width=1\textwidth]{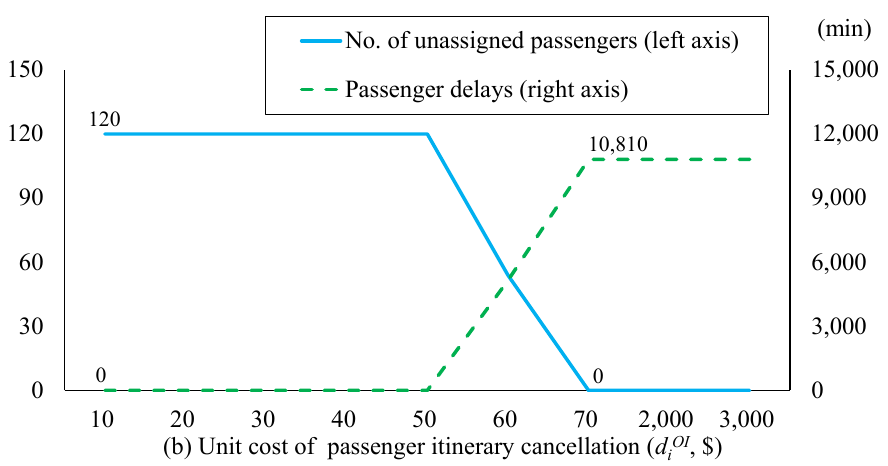}
   \end{minipage}
   \caption{Sensitivity analysis for the unit cost of flight delay and the unit cost of passenger cancellation\label{unit flight delay cost and unit passenger cancellation cost}}
\end{figure}

\noindent
Canceling and delaying are two regular recovery options in both aircraft rerouting and passenger reallocation. There is a trade-off between these two recovery options. If decision-makers prefer lower flight delays when recovering disrupted aircraft routes, a higher unit cost of delaying flights will be set. Consequently, some intractable flights may be canceled to achieve low flight delays. Conversely, if decision-makers wish to avoid flight cancellations, a higher unit cost of canceling flights will be set, and more flights may be delayed in the recovery plan. Clearly, the trade-off between canceling and delaying also exists in passenger reallocation. To explore this trade-off, Figures \ref{unit flight delay cost and unit passenger cancellation cost}(a) and \ref{unit flight delay cost and unit passenger cancellation cost}(b) show the performance when varying the unit costs of delaying flights $(c_{f_{uv}})$ and canceling passenger itineraries $(d_{i}^{OI})$, respectively.

As shown in Figure \ref{unit flight delay cost and unit passenger cancellation cost}(a), as the unit cost of delaying flights increases, the total flight delay time gradually decreases, as expected. Meanwhile, more flights are canceled to reduce the high recovery cost incurred by flight delays. On the other hand, as shown in Figure \ref{unit flight delay cost and unit passenger cancellation cost}(b), fewer passenger itineraries are canceled when the unit cost of canceling itineraries increases. Concurrently, the total passenger delay time increases to preserve as many passenger itineraries as possible.

\subsubsection{Sensitivity Analysis for the Unit Cost of Itinerary Change.}
\begin{figure}[!htbp]
\centering
\includegraphics[width=0.5\textwidth]{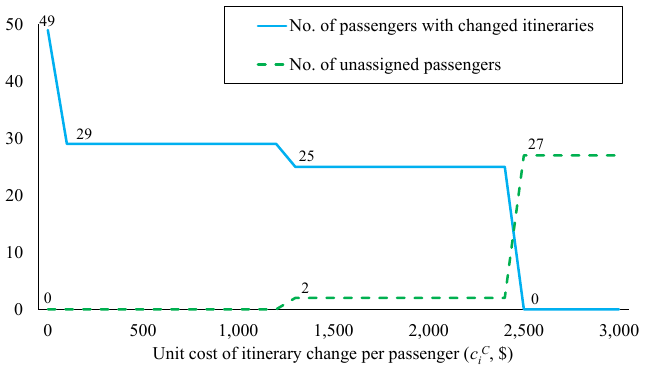}
\caption{Sensitivity analysis for the unit cost of itinerary change \label{unit itinerary change cost}}
\end{figure}

\noindent
Changing passenger itineraries is common in passenger reallocation. However, passengers generally dislike having their planned itineraries changed unless doing so can mitigate the disturbances from disruptions (\citeauthor{arikan2017}, \citeyear{arikan2017}). To examine this situation, as illustrated in Figure \ref{unit itinerary change cost}, we conduct a sensitivity analysis on the unit cost of changing itineraries $(c_{i}^{C})$ to show its effect on recovery performance.

We obtain several observations from Figure \ref{unit itinerary change cost}. First, as the unit cost of changing passenger itineraries increases, the number of passengers with changed itineraries decreases, as expected. Meanwhile, some passengers become unassigned to reduce the recovery cost caused by changing itineraries. Second, when the unit cost of changing itineraries exceeds the unit cost of canceling passenger itineraries (\$2,500 in our experiments), each passenger either follows their planned itinerary or is unassigned in the recovery plan.

\section{Conclusion and future research}\label{Section Conclusion and future research}
In this paper, we introduce a new sparse-dense flight copy approach and propose an interactive mechanism for the integrated flight, aircraft, and passenger recovery problem. The interactive mechanism improves solution efficiency by dividing the integrated recovery problem into two subproblems with much smaller scales and solving them iteratively. To operationalize this mechanism, we develop a path-based integrated flight, aircraft, and passenger recovery model (IFAPRM) and present a customized Benders decomposition (CBD) method to solve it. To further enhance the performance of the CBD method, we propose several acceleration techniques by exploiting the structural properties of the integrated model. Computational experiments demonstrate that the proposed methodology can solve the integrated recovery problem within reasonable runtime, and the acceleration techniques significantly improve the efficiency of the CBD method. The experimental results also validate that the interactive mechanism outperforms the conventional flight copy approach. This interactive mechanism, combined with the mathematical model, the customized Benders decomposition method, and the acceleration techniques, forms a general-purpose methodology that is widely applicable to integrated rescheduling under complex operational disruptions.

This work can be extended in the following directions. First, to enhance the practicality of the proposed interactive mechanism, we will apply it to address other integrated recovery problems. Second, this work does not distinguish between passenger cabin classes to simplify the integrated recovery problem. Future research could consider different cabin classes in passenger reallocation to obtain more economically efficient recovery plans. Third, this work formulates cruise speed control in a discrete manner to alleviate solution difficulty. A continuous formulation could be explored in future work to further improve recovery performance.

% Acknowledgments here
%\ACKNOWLEDGMENT{
%The authors thanks XXX. This work is supported by the National Science Foundation of China (Grant XXX).
%}

% Appendix here
% Options are (1) APPENDIX (with or without general title) or
%             (2) APPENDICES (if it has more than one unrelated sections)
% Outcomment the appropriate case if necessary
%
% \begin{APPENDIX}{<Title of the Appendix>}
% \end{APPENDIX}
%
%   or
%
% \begin{APPENDICES}
% \section{<Title of Section A>}
% \section{<Title of Section B>}
% etc
% \end{APPENDICES}

% References here (outcomment the appropriate case)

% CASE 1: BiBTeX used to constantly update the references
%   (while the paper is being written).
%\bibliographystyle{informs2014trsc} % outcomment this and next line in Case 1
%\bibliography{<your bib file(s)>} % if more than one, comma separated

% CASE 2: BiBTeX used to generate mypaper.bbl (to be further fine tuned)
%\input{mypaper.bbl} % outcomment this line in Case 2

%\bibliographystyle{informs2014trsc}
%\bibliography{ref}

\setlength{\bibsep}{0em}

\newpage
\begin{APPENDICES}
\section{The Proofs of Properties}
\label{appendix proof}
\subsection{Proof of Proposition \ref{BSP2feasibility}}
\noindent
When the BSP is feasible, the BSP2 must be feasible since all the constraints in the BSP2 are also included in the BSP.\par

Conversely, suppose $\widetilde{x}_{f_{uv}}$ is a feasible solution of the BSP2, then we can construct the following solution
\begin{eqnarray*}(x_{f_{uv}}, z_i^{OI}, w_p^{i}) := (\widetilde{x}_{f_{uv}}, n_{i}^{OI}, \mathbf{0}),
\end{eqnarray*}
which is feasible for the BSP.\par

Combining the above two arguments, we conclude the fact that the BSP and the BSP2 have the same feasibility.  \hfill\Halmos

\subsection{Proof of Proposition \ref{Proposition Scale Management 1}}
\noindent
Let $n_{r}^{\bar{R}}$ be the number of flights that are assigned to aircraft $r \in \bar{R}$. It is clear that $n^{\bar{F}} = \sum\limits_{r \in \bar{R}}n_{r}^{\bar{R}}$ holds. Then $n^{Con}$ can be calculated by summing the number of connections in each route, i.e.,
\begin{eqnarray*}
n^{Con} = \sum\limits_{r \in \bar{R}}(n_{r}^{\bar{R}}-1) = \sum\limits_{r \in \bar{R}}n_{r}^{\bar{R}}- \sum\limits_{r \in \bar{R}}1 = n^{\bar{F}}-n^{\bar{R}}.
\end{eqnarray*}  \hfill\Halmos

\vspace{-1em}
\subsection{Proof of Proposition \ref{Proposition Scale Management 2}}
\noindent
The terms $\sum\limits_{f \in F}(1-\bar{z}_f^F)$ and $\sum\limits_{r\in R}\sum\limits_{p\in {P^r}} \bar{y}_p^r$ represent the number of non-canceled flights and the number of selected aircraft, i.e., $n^{\bar{F}}$ and $n^{\bar{R}}$ in Proposition \ref{Proposition Scale Management 1}. Therefore, the number of redundant constraints in (\eqref{FRPRMcon2} is obtained by subtracting the number of selected connections from $|E^{R}|$, i.e.,
\begin{eqnarray*}
|E^{R}|-\left(\sum_{f \in F}(1-\bar{z}_f^F) - \sum_{r\in R}\sum_{p\in {P^r}} \bar{y}_p^r\right).
\end{eqnarray*}  \hfill\Halmos

\vspace{-1em}
\subsection{Proof of Proposition \ref{Proposition Stronger VI}}
\noindent
The proof of the stronger LL cut \eqref{stronger LL cut} is similar to that of the stronger NG cut \eqref{stronger NG cut}, so we only prove the result for the stronger NG cut \eqref{stronger NG cut} as follows.\par

First, we prove that the cut \eqref{stronger NG cut} is valid. We introduce the following two sets related to the RxBMP solution $(\bar{y}_p^r, \bar{z}_f^F)$:
\begin{footnotesize}
\begin{eqnarray*}
&\hspace{-0.2cm}\Omega_1 \hspace{-0.2cm}&= \left\{(y_p^r, z_f^F) \in \{0, 1\}^{\sum\limits_{r\in R}|P^{r}|+|F|}  \Bigg{|} \sum_{(r,p)\in I_{1}(\bar{y}_p^r)}(1-{y_p^r}) + \sum_{(r,p)\in I_{0}(\bar{y}_p^r)}{y_p^r} + \sum_{f\in I_{1}(\bar{z}_f^F)}(1-z_f^F) + \sum_{f\in I_{0}(\bar{z}_f^F)}z_f^F \geq 1\right\},\\
&\hspace{-0.2cm}\Omega_2 \hspace{-0.2cm}&= \left\{(y_p^r, z_f^F) \in \{0, 1\}^{\sum\limits_{r\in R}|P^{r}|+|F|} \Bigg{|} \sum_{(r,p)\in I_{1}(\bar{y}_p^r)}(1-{y_p^r}) + \sum_{(r,p)\in I_{0}(\bar{y}_p^r)}{y_p^r} \geq 1\right\}.
\end{eqnarray*}
\end{footnotesize}\par

Since the NG cut \eqref{NG cut} is valid, then the cut \eqref{stronger NG cut} is also valid if we can prove that $\Omega_1 = \Omega_2$.\par

For any element $(\widetilde{y}_p^r, \widetilde{z}_f^F) \in \Omega_1$, one of the following two cases occurs. Case $(A)$: $\widetilde{z}_f^F$ is equal to $\bar{z}_f^F$ for all $f$. In this case, the term $\sum\limits_{f\in I_{1}(\bar{z}_f^F)}(1-\widetilde{z}_f^F) + \sum\limits_{f\in I_{0}(\bar{z}_f^F)}\widetilde{z}_f^F$ is equal to zero, then the following inequality holds,
\begin{eqnarray*}
\sum\limits_{(r,p)\in I_{1}(\bar{y}_p^r)}(1-{\widetilde{y}_p^r}) + \sum\limits_{(r,p)\in I_{0}(\bar{y}_p^r)}{\widetilde{y}_p^r} \geq 1.
\end{eqnarray*}
Hence, the element $(\widetilde{y}_p^r, \widetilde{z}_f^F)$ belongs to $\Omega_2$. Case $(B)$: $\widetilde{z}_f^F$ is not equal to $\bar{z}_f^F$ for some $f$. In this case, the aircraft rerouting decision in $\widetilde{y}_p^r$ must be different from the decision in $\bar{y}_p^r$, which implies that the element $(\widetilde{y}_p^r, \widetilde{z}_f^F)$ belongs to $\Omega_2$.\par

For any element $(\widetilde{y}_p^r, \widetilde{z}_f^F) \in \Omega_2$, since the term $\sum\limits_{f\in I_{1}(\bar{z}_f^F)}(1-\widetilde{z}_f^F) + \sum\limits_{f\in I_{0}(\bar{z}_f^F)}\widetilde{z}_f^F$ is always nonnegative, it is clear that the element $(\widetilde{y}_p^r, \widetilde{z}_f^F)$ belongs to $\Omega_1$. Hence, $\Omega_1 = \Omega_2$, and it follows that the cut \eqref{stronger NG cut} is valid.\par

\vspace*{0.5cm}
Second, we show that the cut \eqref{stronger NG cut} is stronger than the NG cut \eqref{NG cut}. Recall that $(y_p^r, z_f^F)$ are binary in the Benders master problem. Then, when $(y_p^r, z_f^F)$ are relaxed to $[0, 1]^{\sum\limits_{r\in R}|P^{r}|+|F|}$, the following inequality holds
\begin{footnotesize}
\begin{eqnarray*}
&&\hspace{-0.5cm}\sum_{(r,p)\in I_{1}(\bar{y}_p^r)}(1-{y_p^r}) + \sum_{(r,p)\in I_{0}(\bar{y}_p^r)}{y_p^r}\leq\sum_{(r,p)\in I_{1}(\bar{y}_p^r)}(1-{y_p^r}) + \sum_{(r,p)\in I_{0}(\bar{y}_p^r)}{y_p^r} + \sum_{f\in I_{1}(\bar{z}_f^F)}(1-z_f^F) + \sum_{f\in I_{0}(\bar{z}_f^F)}z_f^F,
\end{eqnarray*}
\end{footnotesize}
\noindent
since the term $\sum\limits_{f\in I_{1}(\bar{z}_f^F)}(1-z_f^F) + \sum\limits_{f\in I_{0}(\bar{z}_f^F)}z_f^F$ is always nonnegative. Then it follows that the valid inequality \eqref{stronger NG cut} is stronger than the NG cut \eqref{NG cut}.  \hfill\Halmos

\vspace{2em}
\section{Constraint Programming Models}\label{Section Constraint programming models}
\noindent
In this section, we provide the constraint programming (CP) models of passenger itineraries and aircraft routes. Before giving the specific formulations, the involved parameters and functions are presented as follows.\par

\begin{center}
\renewcommand{\baselinestretch}{1.15}
\begin{small}
    \vspace*{-0.5cm}
	\begin{longtable}{p{0.25\textwidth}p{0.63\textwidth}}
        \caption{Parameters used in the CP models}\label{Parameters used in the CP models}
        \vspace*{-0.2cm}
		\endfirsthead
			
        \midrule
		\endhead
			
        \midrule
        \multicolumn{2}{l}{\footnotesize Continued on next page}
		\endfoot
			
		\hline
	\endlastfoot

    \toprule
    $SDepAP_{f_{u}}/SArrAP_{f_{u}}$ & the departure/arrival airport of sparse flight copy ${f_{u}}$\\
    $SDepT_{f_{u}}/SArrT_{f_{u}}$ & the departure/arrival time of sparse flight copy ${f_{u}}$\\
    $DDepAP_{f_{uv}}/DArrAP_{f_{uv}}$ & the departure/arrival airport of dense flight copy ${f_{uv}}$\\
    $DDepT_{f_{uv}}/DArrT_{f_{uv}}$ & the departure/arrival time of dense flight copy ${f_{uv}}$\\
    $StartAP_{e}/EndAP_{e}$ & the start/end airport of entity (original itinerary or aircraft) ${e}$\\
    $TurnT_{r}$ & the minimum turn time of aircraft $r$\\
    $SMT_{r}/DMT_{r}$ & the index of the sparse/dense flight copy corresponding to the maintenance task of aircraft $r$\\
    $SFuelC_{r,f_{u}}$ & the change of the used fuel if sparse flight copy ${f_{u}}$ is performed by aircraft $r$\\
    $SSwap_{r,f_{u}}$ &1, if flight $f$ is an unplanned flight for aircraft $r$;~0, otherwise\\
    $StartT_{i}/EndT_{i}$ & the planned start/end time of original itinerary ${i}$\\
    $DSwap_{i,f_{uv}}$ &1, if flight $f$ is the unplanned flight for original itinerary $i$;~0, otherwise\\
    $SDual_{f_{u}}$ & the dual value related to sparse flight copy ${f_{u}}$\\
    $DDual_{f_{uv}}$ & the dual value related to dense flight copy ${f_{uv}}$\\
    $SPairDual_{(f_{1u}, f_{2u})}$ & the dual value related to the pair of sparse flight copies ${(f_{1u}, f_{2u})}$\\
    $RDual_{r}$ & the dual value related to aircraft $r$\\
    $OIDual_{i}$ & the dual value related to original itinerary $i$\\
    $SDepAP/SArrAP$ & the array of $SDepAP_{f_{u}}/SArrAP_{f_{u}}$\\
    $SDepT/SArrT$ & the array of $SDepT_{f_{u}}/SArrT_{f_{u}}$\\
    $DDepAP/DArrAP$ & the array of $DDepAP_{f_{uv}}/DArrAP_{f_{uv}}$\\
    $DDepT/DArrT$ & the array of $DDepT_{f_{uv}}/DArrT_{f_{uv}}$\\
    $SFuelC_{r}/SSwap_{r}$ & the array of $SFuelC_{r,f_{u}}/SSwap_{r,f_{u}}$ for aircraft $r$\\
    $DSwap_{i}$ & the array of $DSwap_{i,f_{uv}}$ for original itinerary $i$\\
    $SDual/DDual$ & the array of $SDual_{f_{u}}/DDual_{f_{uv}}$\\
    $SPairDual$ & the array of $SPairDual_{(f_{1u}, f_{2u})}$\\
    $c_{\textrm{swap}}$ & the unit cost of swapping tail assignments\\
    $c_{\textrm{fuel}}$ & the unit fuel cost\\
    $\kappa$ & the CO$_{2}$ emission constant\\
    $c_{\textrm{CO$_2$}}$ & the unit cost of carbon emission\\
    $FC^{D}$ & the set of all dense flight copies\\
    $ConT$ & the minimum connection time for passengers\\
    \bottomrule
	\end{longtable}
    \vspace*{-0.5cm}
\end{small}
\end{center}

The dual values mentioned in Table \ref{Parameters used in the CP models} can be obtained according to Equation \eqref{reduced cost itinerary} and Equation \eqref{reduced cost route}. Then, we provide the definitions of functions $presenceOf(\cdot),~count(\cdot, \cdot),$ and $element(\cdot, \cdot)$ in constraint programming. Suppose the decision variable $p=(p_{1}, p_{2}, ..., p_{\eta})$ is a vector of the indexes of dense/sparse flight copies indicating the flight sequence of the associated passenger itineraries/aircraft routes, where $\eta$ denotes the maximum allowable number of legs. The function $presenceOf(p_{k})$ indicates whether the element $p_{k}$ is present. The function $count(p, \alpha)$ returns the number of the elements that are equal to the value $\alpha$ in the vector $p$. Suppose the value of $p_{k}$ is equal to $\alpha$, then the function $element(X, p_{k})$ returns the $\alpha$th element of the array $X$. Accordingly, suppose the values of $p_{k}$ and $p_{k+1}$ are equal to $\alpha$ and $\beta$ respectively, then the function $element(Y, (p_{k}, p_{k+1}))$ returns the ($\alpha$th, $\beta$th) element of the two-dimensional array $Y$.\par

\subsection{The Constraint Programming Model of Passenger Itineraries}\label{Section The constraint programming model of passenger itineraries}
\noindent
We now provide the CP model of the passenger itineraries for original itinerary $i$. In this model, the decision variable $p=(p_{1}, p_{2}, ..., p_{\eta}) \in \{0, 1, ..., |FC^{D}|-1\}^{\eta_{i}}$ is a vector of the indexes of dense flight copies, where $\eta_{i}$ denotes the maximum allowable number of legs for the passengers in original itinerary $i$. The specific constraints are presented as follows.\par

\vspace{-10pt}
\begin{eqnarray}
&\hspace{-10cm}\textbf{\small Constraints:}\notag\\
&presenceOf(p_{k+1}) \leq presenceOf(p_{k}),~~\forall k\in \{1, ..., \eta_{i}-1\},\label{CP2 con1}\\
&presenceOf(p_{k+1}) == (p_{k+1}!=p_{k}),~~\forall k\in \{1, ..., \eta_{i}-1\},\label{CP2 con2}\\
&element(DDepAP, p_1) == StartAP_{i},\label{CP2 con3}\\
&element(DArrAP, p_{\eta_{i}}) == EndAP_{i},\label{CP2 con4}\\
&\left[element(DArrT, p_{k}) + ConT \leq element(DDepT, p_{k+1})\right]\notag\\
&~~~~~~\bigvee~\left[presenceOf(p_{k+1}) == 0\right],~~\forall k\in \{1, ..., \eta_{i}-1\}, \label{CP2 con6}\\
&\left[element(DArrS, p_{k}) == element(DDepS, p_{k+1})\right]\notag\\
&~~~~~~\bigvee~\left[presenceOf(p_{k+1}) == 0\right],~~\forall k\in \{1, ..., \eta_{i}-1\}, \label{CP2 con7}\\
&count(p, DMT_{r}) == 0,~~\forall r\in R, \label{CP2 con5}\\
&element(DDepT, p_1) \geq StartT_{i},,\label{CP2 con8}\\
&DT^{i}_{p} = \max\left\{0,~element(DArrT, p_{\eta}) - EndT_{i} \right\},\label{CP2 con9}\\
&n^{C}_{i,p} = \sum\limits_{k=1}^{\eta} element(DSwap_{i}, p_{k})\cdot presenceOf(p_{k}),\label{CP2 con10}\\
&\hspace*{-2cm}\bar{c}_p^{i} = c_{i}^{DT}DT^{i}_{p} + c_{i}^{C}n^{C}_{i,p} - OIDual_{i} \notag\\
&\hspace*{2cm}+ \sum\limits_{k=1}^{\eta} element(DDual, p_{k})\cdot presenceOf(p_{k}).\label{CP2 con11}
\end{eqnarray}

Constraints \eqref{CP2 con1} ensure that $p_{k}$ must be present if $p_{k+1}$ is present. With the aid of Constraints \eqref{CP2 con1}, the value of $presenceOf(p_{k+1})$ also indicates whether the connection of sparse flight copies $(p_{k}, p_{k+1})$ is present. The value of $presenceOf(p_{k+1})$ is equal to one if and only if the connection $(p_{k}, p_{k+1})$ is present in the vector $p$. Constraints \eqref{CP2 con2} mean that the value of $p_{k+1}$ is different from that of $p_{k}$ if $p_{k+1}$ is present; otherwise, the value of $p_{k+1}$ is the same as that of $p_{k}$. Constraints \eqref{CP2 con3} and \eqref{CP2 con4} restrict the start airport and the end airport for original itinerary $i$, respectively. Combined with Constraints \eqref{CP2 con1}, Constraints \eqref{CP2 con6} and \eqref{CP2 con7} make sure that the minimum connection time and the space match are satisfied for each present flight connection in the vector $p$. Constraints \eqref{CP2 con5} ensure that passenger itineraries does not contain the dense flight copies corresponding maintenance tasks. Constraint \eqref{CP2 con8} ensures that the start time of the itinerary should be later than the planned start time of the associated original itinerary. Equations \eqref{CP2 con9} and \eqref{CP2 con10} calculate the delay time and the number of the unplanned flights for the passenger itinerary corresponding to the vector $p$, respectively. Equation \eqref{CP2 con11} calculates the reduced cost of the associated passenger itinerary.

\subsection{The Constraint Programming Model of Aircraft Routes}\label{Section The constraint programming model of aircraft routes}
\noindent
In this subsection, we provide the CP model of the routes for aircraft $r$. In this model, the decision variable $p=(p_{1}, p_{2}, ..., p_{\eta_{r}}) \in \{lb_r, lb_r + 1, ..., ub_r\}^{\eta_{r}}$ is a vector of the indexes of sparse flight copies that can be performed by aircraft $r$, where $\eta_{r}$ denotes the maximum allowable number of legs for aircraft $r$. The specific constraints of the routes for aircraft $r$ are presented as follows.\par

\vspace{-10pt}
\begin{eqnarray}
&\hspace{-10cm}\textbf{\small Constraints:}\notag\\
&presenceOf(p_{k+1}) \leq presenceOf(p_{k}),~~\forall k\in \{1, ..., \eta_{r}-1\},\label{CP1 con1}\\
&presenceOf(p_{k+1}) == (p_{k+1}!=p_{k}),~~\forall k\in \{1, ..., \eta_{r}-1\},\label{CP1 con2}\\
&element(SDepAP, p_1) == StartAP_{r},\label{CP1 con3}\\
&element(SArrAP, p_{\eta_{r}}) == EndAP_{r},\label{CP1 con4}\\
&\left[element(SArrT, p_{k}) + TurnT_{r} \leq element(SDepT, p_{k+1})\right]\notag\\
&~~~~~~\bigvee~\left[presenceOf(p_{k+1}) == 0\right],~~\forall k\in \{1, ..., \eta_{r}-1\}, \label{CP1 con7}\\
&\left[element(SArrS, p_{k}) == element(SDepS, p_{k+1})\right]\notag\\
&~~~~~~\bigvee~\left[presenceOf(p_{k+1}) == 0\right],~~\forall k\in \{1, ..., \eta_{r}-1\}, \label{CP1 con8}\\
&count(p, SMT_{r}) \geq 1,\label{CP1 con5}\\
&count(p, SMT_{\bar{r}}) == 0,~~\forall \bar{r}\in R\setminus\{r\}, \label{CP1 con6}\\
&\textrm{SwapN} = \sum\limits_{k=1}^{\eta_{r}} element(SSwap_{r}, p_{k})\cdot presenceOf(p_{k}), \label{CP1 con9}\\
&\textrm{FuelC} = \sum\limits_{k=1}^{\eta_{r}} element(SFuelC_{r}, p_{k}) \cdot presenceOf(p_{k}), \label{CP1 con10}\\
&c_p^r =  c_{\textrm{swap}} \cdot \textrm{SwapN} + (c_{\textrm{fuel}} + \kappa\cdot c_{\textrm{CO$_2$}}) \cdot \textrm{FuelC}, \label{CP1 con11}\\
&\hspace*{-0.5cm}\bar{c}_p^r = c_p^r + RDual_r + \sum\limits_{k=1}^{\eta_{r}} element(SDual, p_{k})\cdot presenceOf(p_{k}) \notag\\
&\hspace*{0.5cm}+ \sum\limits_{k=1}^{\eta_{r}-1} element(SPairDual, (p_{k},p_{k+1}))\cdot presenceOf(p_{k+1}). \label{CP1 con12}
\end{eqnarray}

Constraints $\eqref{CP1 con1}-\eqref{CP1 con8}$ are similar to Constraints $\eqref{CP2 con1}-\eqref{CP2 con7}$ in the CP model of passenger itineraries. Constraint \eqref{CP1 con5} guarantees the maintenance requirement is satisfied for aircraft $r$. Constraints \eqref{CP1 con6} make sure that the route of aircraft $r$ can not contain the maintenance tasks of other aircraft. Equations \eqref{CP1 con9} and \eqref{CP1 con10} calculate the number of the unplanned flights and the change of the used fuel for the aircraft route corresponding to the vector $p$, respectively. The cost and the reduced cost of the associated aircraft route are calculated by Equation \eqref{CP1 con11} and \eqref{CP1 con12}.

\end{APPENDICES}

\end{document}